\documentclass[a4paper,11pt]{article}

\usepackage{subfiles}

\usepackage[english]{babel}
\usepackage[utf8]{inputenc}
\usepackage[T1]{fontenc}
\usepackage{lmodern}
\usepackage{amsmath}
\usepackage{amssymb}

\usepackage{amsthm}
\usepackage{booktabs}
\usepackage{subcaption}
\usepackage{xcolor}
\usepackage{graphicx}
\usepackage{bm}
\usepackage{enumitem}

\usepackage[top=25mm, bottom=25mm, left=25mm, right=25mm]{geometry}
\usepackage{setspace}
\usepackage[round]{natbib}

\usepackage[ruled]{algorithm2e}
\usepackage{sectsty}
\usepackage{multirow}
\usepackage{titling}
\usepackage{todonotes}
\usepackage{scrextend}

\usepackage{mathtools}

\DeclarePairedDelimiter\floor{\lfloor}{\rfloor}

\DeclareCaptionFormat{part1}{#1 (part 1)#2#3\par}
\DeclareCaptionFormat{part2}{#1 (part 2)#2#3\par}

\usepackage{hyperref}
\usepackage{cleveref}
\hypersetup{
hidelinks
}

\newcommand{\orcid}[1]{\href{https://orcid.org/#1}{\includegraphics[scale=0.12]{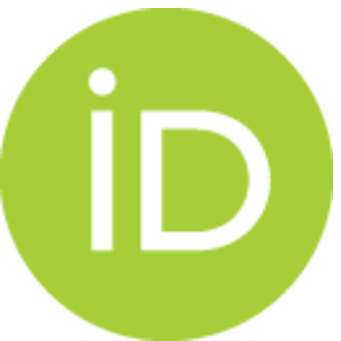}}}

\DeclareMathOperator*{\argmax}{arg\,max}
\DeclareMathOperator*{\argmin}{arg\,min}

\newlist{CC}{enumerate}{1}
\setlist[CC]{label=Case \arabic*)}

\newlist{STP}{enumerate}{1}
\setlist[STP]{label=Step \arabic*)}

\usepackage{thmtools}
\theoremstyle{plain}
\newtheorem{lemma}{Lemma}

\theoremstyle{definition}
\newlist{asslist}{enumerate}{1}
\setlist[asslist]{label=(\roman{asslisti}), ref=\theassumption(\roman{asslisti}),noitemsep}
\declaretheorem[
name=Assumption,
Refname={Assumption,Assumptions}]{assumption}
\Crefname{asslisti}{Assumption}{Assumptions}
\addtotheorempostheadhook[assumption]{\crefalias{asslisti}{assumption}}

\newcommand \linedabstractkw[2]{% Default width = 0.9
  \renewcommand\maketitlehookd{%
    \mbox{}\medskip\par
    \centering
    \hrule\medskip
    \begin{minipage}{0.9\textwidth}
    #1\\

    \textit{Keywords: }#2
    \end{minipage}\medskip\hrule\medskip
    }
}

\graphicspath{{img/}{../img/}}

\title{
Resource allocation problems with expensive function evaluations
}
\author{S.C.M. ten Eikelder\thanks{Department of Econometrics and Operations Research, Tilburg University, The Netherlands. \url{https://orcid.org/0000-0001-7883-8506}}~~\orcid{0000-0001-7883-8506} \and J.H.M. van Amerongen\thanks{Independent researcher, Leiden, The Netherlands. \url{https://orcid.org/0000-0002-1706-8779}}~~\orcid{0000-0002-1706-8779}
}

\date{May 3, 2022}

\begin{document}
%%%%%%%%%%%%%%%%

\linedabstractkw{The resource allocation problem is among the classical problems in operations research, and has been studied extensively for decades. However, current solution approaches are not able to efficiently handle problems with expensive function evaluations, which can occur in a variety of applications. We study the integer resource allocation problem with expensive function evaluations, for both convex and non-convex separable cost functions. We present several solution methods, both heuristics and exact methods, that aim to limit the number of function evaluations. The methods are compared in numerical experiments using both randomly generated instances and instances from two resource allocation problems occurring in radiation therapy planning. Results show that the presented solution methods compare favorably against existing derivative free optimization solvers.}{Nonlinear programming, resource allocation problem, expensive function evaluations, black-box optimization, radiation therapy}

\maketitle

\section{Introduction}\label{sec: introduction}
Resource allocation problems are among the classical problems in operations research, with the earliest investigations in the 1950s \citep{Koopman53}. In a generic resource allocation problem, a decision maker has a fixed amount of resources, and the goal is to divide these over a set of players, tasks or projects such that a cost function is minimized. Many variations of this problem have been studied in literature, with either continuous or integer variables, objective functions that are separable or non-separable and convex or non-convex. The unifying characteristic for most of these formulations is a single (linear or non-linear) constraint on the total amount of resource to be allocated, aside of variable bounds. Thus, resource allocation problems can be seen as a special case of nonlinear (integer) programming\footnote{In literature, resource allocation problems are also referred to as nonlinear knapsack problems.} \citep{Hochbaum07}. Resource allocation problems are amongst others encountered in production and inventory management, economics, finance, allocation of computer resources, and telecommunications. We refer to \citet{Bretthauer02}, \citet{Patriksson08}, \citet{Katoh13}, and \citet{Patriksson15} for reviews of problem formulations, algorithms and applications. Recent applications also include vaccine allocation in epidemiology \citep{Duijzer18} and decentralized energy management \citep{SchootUiterkamp21}. The latter also provides an overview of algorithms and complexity results for various problem formulations.

In current solution approaches, it is generally assumed that the cost associated with a particular resource allocation is easily computed. However, determining the cost (or value) of allocating a certain number of resources to an individual project or player may be non-trivial in practical applications; this may be expensive, either financially or time-wise. In such cases, the resource allocation problem is said to have \emph{expensive function evaluations}. If funds or time are limited, the problem cannot readily be solved using existing solution approaches. The purpose of this paper is to introduce the resource allocation problem with expensive function evaluations, and present solution methods. 

\subsection{Contributions}
We present methods to find optimal or near-optimal solutions to integer resource allocation problems, while limiting the number of function evaluations. To the best of our knowledge, the (integer) resource allocation problem with expensive function evaluations has not been studied in literature. We consider problem formulations with both convex and non-convex non-increasing cost functions. 

We propose two novel solution methods. The first, the 1-Opt method, starts with a feasible allocation and subsequently evaluates new points where a single move (1-opt step) potentially leads to improved allocations. This method is  exact for convex cost functions, and a heuristic for non-convex cost functions. The second proposed method, the sandwich method, also guarantees a globally optimal solution for non-convex cost functions. It sandwiches the cost function between an upper and lower bound, and evaluates those points that are expected to furthest reduce this gap. Both methods are compared to several benchmark methods, including the NOMAD solver \citep{NOMAD}.

The performance of the methods is compared on various randomly generated instances, both with convex and non-convex cost functions. We also investigate the influence of early termination of the methods, i.e., what intermediate solution and objective value bounds can be obtained when the method is manually terminated by the decision-maker. Furthermore, we consider two applications from radiation therapy planning with `near-convex' cost functions. To summarize, our contributions are the following:
\begin{itemize}
\item We formulate the resource allocation problem with expensive function evaluations, and present heuristic and exact solution methods that aim to limit the number of function evaluations.
\item We present numerical experiments on randomly generated instances with convex and non-convex cost functions, and instances from two radiation therapy applications. Whereas the sandwich method performs best on both the convex and non-convex randomly generated instances, the 1-Opt method performs best on the near-convex radiation therapy instances. Both methods consistently outperform the benchmark methods.
\end{itemize}

\subsection{Applications}
The resource allocation problem with expensive function evaluations has several areas of application. Two examples are capital investment and radiation therapy planning.

In making strategic capital investment decisions, decision makers often have to allocate a certain research or marketing budget over a number of projects. For an individual project, finding the most efficient way to spend this money may require costly simulations, expensive market potential research or consultancy costs. Thus, one aims to limit the number of such initial studies during the decision-making process.
For example, in the Netherlands a quantitative approach to flood protection has been adopted \citep{Eijgenraam14}. For each dike-ring area, 53 in total, one can determine how to efficiently improve flood protection standards given a certain budget, but these studies are time-expensive.

In radiation therapy planning for cancer treatments, we have encountered resource allocation problems with expensive function evaluations in two situations. These are studied in the numerical experiments in \Cref{sec: RT}. 
\begin{enumerate} 
\item Proton therapy is an expensive and in many countries scarce radiation therapy modality, and the available treatment slots should be allocated to those patients who are expected to benefit the most. Determining the optimal radiation therapy treatment plan for a cancer patient given a certain amount of proton slots is computationally expensive, and these allocation decisions need to be made on a weekly basis for potentially large patient populations.
\item During the planning process of volumetric modulated arc therapy (VMAT), a particular radiation therapy delivery method, the gantry (i.e., the treatment device) rotates around the patient while continuously irradiating the patient. One typically tries to limit treatment delivery time, which in some approaches translates to an upper bound on the total available treatment time. This time needs to be efficiently allocated to different segments of the entire 360 degree arc. The optimal delivery plan for each arc segments depends on the allocated treatment time, and is computationally expensive to compute.
\end{enumerate}

\subsection{Problem formulation}
The resource allocation problem aims to allocate a finite set of identical items over a set of players, indexed $i=1,\dotsc,n$, such that a cost function is minimized. Typically, this cost function is separable. Let an allocation be represented by decision variable $\bm{x} \in \mathbb{N}_+^n$, and let $B\in \mathbb{N}_+$ denote the total allocation budget, $\bm{b}\in \mathbb{N}_+^n$ the individual budgets and $f_i : \{0,\dotsc,b_i\} \mapsto \mathbb{R}$ the non-increasing cost function for player $i$. Then the resource allocation problem reads
\begin{subequations}\label{eq: knapsack}
\begin{align}
P(f)~~~~\min_{\bm{x}}~&~ \sum_{i=1}^n f_i(x_i) \\
\text{s.t.}~&~ \sum_{i=1}^n x_i = B \label{eq: knapsack-total-budget}\\
           ~&~ 0 \leq x_i \leq b_i,~x_i\text{ integer},~\forall i=1,\dotsc,n.
\end{align}
\end{subequations}
This problem is referred to as $P(f)$, where a problem instance is represented by the objective function $f : \mathbb{R}^n \mapsto \mathbb{R}$ defined by $f(\bm{x}) = \sum_i f_i(x_i)$.\footnote{Note that a full instance of problem $P(\cdot)$ is specified by the objective function $f$, individual budgets $\bm{b}$ and the total allocation budget $B$. We represent $\bm{b}$ in $f$ via the domains of $f_i$ and omit $B$ for notational convenience.} The optimal objective value of problem $P(f)$ is denoted $z^{\ast}(f)$, the objective value of a feasible solution $\bm{x}$ evaluated on $f$ is denoted $z(\bm{x},f)$.

We study cases in which every individual cost function $f_i$ is deterministic yet unknown a priori (black-box) and for which it is expensive to determine the cost of allocating a certain number of resources to a certain player, which we call a \emph{function evaluation}\footnote{Note that a function evaluation refers to the evaluation of an individual cost function $f_i(\cdot)$ at some value $x_i$, not to the evaluation of the entire cost function $f = \sum_i f_i$ of optimization problem $P(f)$.}.
\begin{assumption}\label{ass: expensive-evaluations}
The expense of evaluating $f_i(x_i)$ for any $x_i$ and $i=1,\dotsc,n$ is much larger than $(\gg)$ the expense of solving $P(\cdot)$ with known function values.
\end{assumption}
\noindent \textbf{Goal:} To find an optimal solution to $P(f)$ while using as few function evaluations of functions $f_i$, $i=1,\dotsc,n$, as possible.\\

Without knowing the cost functions, $P(f)$ is not fully specified and cannot be solved directly. We propose methods that solve a sequence of subproblems that use partial information available on functions $f_i(\cdot)$, $i=1,\dotsc,n$. This information is captured in a set of deterministic, known cost functions $c_i :\{0,\dotsc,b_i\} \mapsto \mathbb{R}$, $i=1,\dotsc,n$ and the associated, fully specified and solvable, resource allocation problem is denoted $P(c)$.

\subsection{Assumptions}
We make several technical assumptions to ensure Problem \ref{eq: knapsack} is both interesting and solvable. We assume that the total allocation budget is restrictive (i.e. not each player $i$ can be allocated their budget $b_i$), so that all items are allocated in an optimal allocation. For ease of exposition it is additionally assumed that the individual budgets do not exceed the total allocation budget.
Furthermore, we assume to know bounds on functions $f_i$. For ease of exposition, we let these bounds be $[0,M]$ for every function, where $M\in\mathbb{R}_+$ is some large number; other bounds may be chosen. These assumptions are summarized below.
\begin{assumption} \label{ass: technical}
It holds that\footnote{
These conditions also apply when $P(\cdot)$ is solved using any other generic cost function $\bm{c}$, where $f_i$ in \mbox{\Cref{ass: bounds}} is replaced by $c_i$.
}
\leavevmode
\begin{asslist} 
\item $\bm{b}^{\top}\bm{e} > B$
\item $B \geq b_i$ for all $i=1,\dotsc,n$
\item $0\leq f_i(x_i) \leq M$ for all $x_i \in \{0,\dotsc,b_i\}$, $i=1,\dotsc,n$. \label{ass: bounds}
\end{asslist}
\end{assumption}

The naive approach is to compute for each player $i$ the values of $f_i(x_i)$ for all feasible $x_i$, and solve $P(f)$ using any standard method. This requires $\bm{b}^{\top}\bm{e}+n$ expensive function evaluations. Any solution method (heuristic, approximation or exact) should have limited function evaluations to be of practical value.

The difficulty of finding the optimal solution to $P(f)$ using few function evaluations depends on what information is available concerning the behavior of the functions $f_i$. We consider two extreme cases. First, we assume that all functions $f_i$, $i=1,\dotsc,n$, are known to be convex and non-increasing in $x_i$. This may for example occur in resource allocation problems in marketing, where projects often exhibit diminishing returns to scale. Afterwards we assume all functions $f_i$, $i=1,\dotsc,n$, are solely known to be non-increasing in $x_i$. In both cases, we assume no other second-order or probabilistic information on the behavior of functions $f_i$ is known. In some applications, it may be known that the cost of computing $f_i(k)$ is not constant, e.g., it might increase in $k$. In such cases, minimizing total evaluation expenses is not equivalent to minimizing the total number of function evaluations. We assume not to have such information available.

\subsection{Literature review} \label{sec: literature_review}
Optimization problems with expensive function evaluations can also be solved using derivative-free optimization (DFO) approaches \citep{Audet17,Larson19}. These methods have been developed for many problem types where derivatives cannot be computed or approximated. This is amongst others the case if the objective function is a black-box, e.g., if it is evaluated via computer simulations. Some DFO methods, e.g., \citet{Brekelmans05}, also try to limit the number of function evaluations. 

The majority of DFO research focuses on optimization problems with solely continuous variables; integer and mixed-integer approaches are more recent. \citet{Ploskas21} provide a review of DFO algorithms and software for mixed-integer problems. Most of the presented algorithms work for problems with only bound constraints, although linear constraints (such as the resource allocation budget constraint) can be included via penalty objective terms. \citet{Larson21} consider derivative-free minimization of a convex function on an integer lattice. They use an underestimator that interpolates between previously evaluated points; the underestimator determines new points to be evaluated, until global optimality is certified. Their approach does not use overestimators and does not exploit separability of the objective function. The sandwich method proposed in the current paper is similar, but alleviates these two limitations.

A drawback of ignoring separability is that the solution space is vastly larger. In problem \eqref{eq: knapsack}, there are $n + \bm{b}^{\top}\bm{e}$ data points $f_i(x_i)$ that may be evaluated. Ignoring separability, one can distinguish $\prod_{i=1}^n (b_i + 1)$ function values $f(\bm{x})$, which grows exponentially in $n$.

In \Cref{sec: numerical-experiments} the open-source DFO solver NOMAD \citep{NOMAD} is used as a benchmark for the newly proposed solution methods. In preliminary numerical experiments, we also tested MISO \citep{Muller16}, BFO \citep{Porcelli17}, and Matlab's genetic algorithm (GA) implementation\footnote{\url{https://mathworks.com/help/gads/ga.html}.}, but these yielded unsatisfactory results. MISO starts with an intial Latin Hypercube experimental design. Because such designs are typically non-collapsing on the $n$-dimensional feasible region, every individual cost function $f_i$ is already evaluated on most of its domain $\{0,\dotsc,b_i\}$. Similarly, Matlab's GA uniformly samples an initial population from the feasible region, which also results in (near)-complete evaluation of all individual cost functions.

The BFO solver, when applied to a pure integer problem, considers an individual variable $x_i$, and attempts to increase or decrease its value until no further improvement is possible, and then moves to the next variable. For the current problem, this means that each variable $x_1, x_2, \dotsc$ is set at value $b_1, b_2, \dotsc$ until the total budget $B$ is reached. No local search step is performed. Preliminary numerical experiments show that this brute-force approach is outperformed by other benchmark methods.

Altogether, because DFO methods are not developed specifically for resource allocation problems, and/or do not use separability of cost functions, they may be outperformed by dedicated solution approaches.

\section{Solution methods} \label{sec: methods}
In this section two solution methods are presented for solving or approximating $P(f)$ using few expensive function evaluations. Their solution quality guarantees are discussed for both convex and non-convex cost functions. First, we define valid lower and upper bounds on the cost functions, given a set of evaluated points, for convex and non-convex cost functions. Subsequently, we describe the 1-Opt method and the sandwich method.

\subsection{Bounds for non-convex cost functions} \label{sec: properties-general}
\begin{figure}[h!]
\begin{subfigure}{0.5\textwidth}
\includegraphics[scale=1]{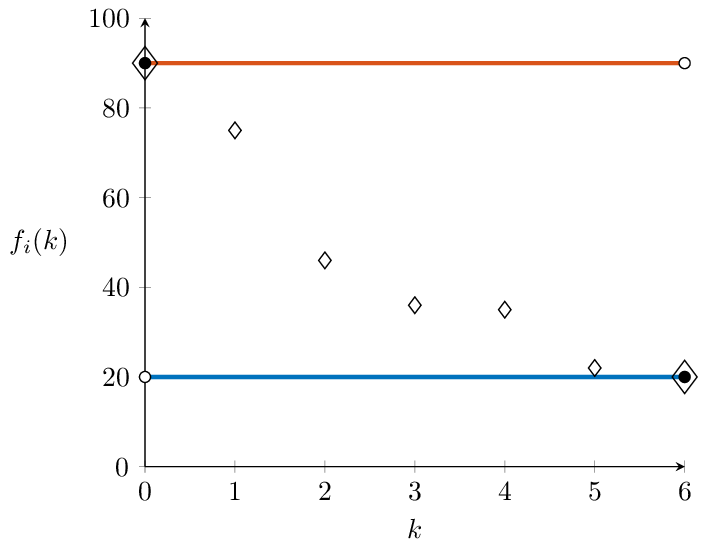}
\caption{\small Two data points. \label{fig: piecewise-general2}}
\end{subfigure}
\begin{subfigure}{0.5\textwidth}
\includegraphics[scale=1]{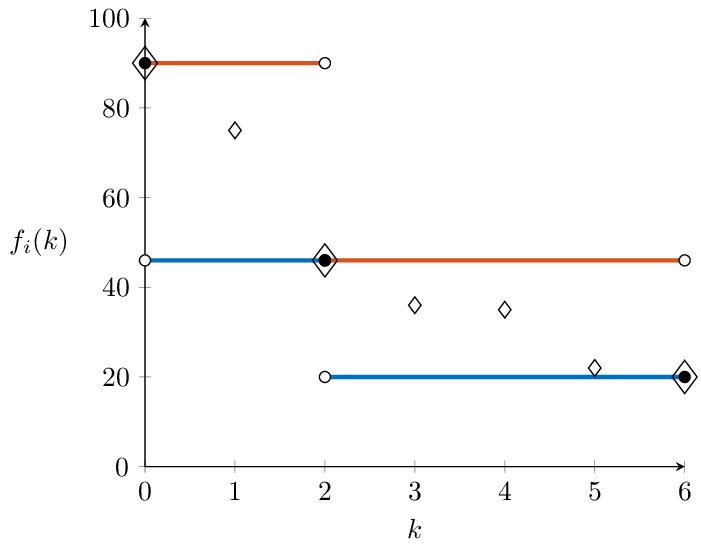}
\caption{\small Three data points.\label{fig: piecewise-general3}}
\end{subfigure}
\caption{\small Non-convex cost functions. Evaluated and non-evaluated points are indicated by large and small diamonds, respectively. If $f_i(k)$ is evaluated for only a few values of $k$, we have an upper (red) and lower (blue) bound for the non-evaluated points. Discontinuity points are indicated by circles.\label{fig: piecewise-general}}
\end{figure}

If for a player $i$ not all points on the cost curve are evaluated, we can obtain a lower and upper bound for each point based on the points that we do have evaluated, with more points yielding tighter bounds. \Cref{fig: piecewise-general} gives an example for a non-convex non-increasing cost function $f_i: \{0,\dotsc,b_i\} \rightarrow [0,M]$, with $b_i=6$. In \Cref{fig: piecewise-general2} two points are evaluated: $f_i(0)$ and $f_i(6)$. The horizontal blue and red line illustrates that the values $f_i(6)$ and $f_i(0)$ are lower and upper bounds for $f_i(k)$ for all $k$, respectively, because $f_i$ is non-increasing. In \Cref{fig: piecewise-general3} additionally $f_i(2)$ is evaluated. The extra evaluated point improves the upper or lower bound for all non-evaluated points $f_i(k)$, $k=0,\dotsc,b_i$. Large gaps between the upper and lower bound on a particular $f_i(k)$ suggest that computing that point yields much information on the true shape of $f_i$. However, it is important to note that this information does not necessarily contribute to solving $P(f)$; some parts of $f_i$ may be irrelevant for the optimal solution.
We proceed by constructing lower and upper bounds. Define the evaluation indicator
\begin{align}
v_{i,k} = 
\begin{cases}
1 & \text{ if } f_i(k) \text{ has been evaluated} \\
0 & \text{otherwise,}
\end{cases}
\end{align}
and let $V$ denote the matrix with elements $v_{i,k}$. For a given evaluation matrix $V$ one can construct lower and upper bound cost functions $l_i, u_i: \{0,\dotsc,b_i\} \mapsto \mathbb{R}$ for each player $i$\footnote{Functions $l_i$ and $u_i$ depend on evaluation matrix $V$ but for notational convenience this is omitted.}. The lower and upper bounds on $f_i(k)$ are given by
\begin{subequations}
\begin{align}
l_{i}(k) &= \max\{0,\max_{p\geq k} f_i(p)v_{i,p}\} \\
u_{i}(k) &= \min\{M,\min_{q\leq k} f_i(q)v_{i,q}\}.
\end{align}
\end{subequations}
If $v_{i,p}=1$ for some $p\geq k$ then the lower bound is always attained at the smallest such $p$ because $f_i(p)$ is decreasing in $p$. For the same reason, if $v_{i,q}=1$ for some $q\leq k$ then the upper bound is always attained at the largest such $q$. For evaluated points $f_i(k)$, the lower and upper bounds coincide. Let $l(\bm{x}) := \sum_i l_i(x_i)$ and $u(\bm{x}) := \sum_i u_i(x_i)$. Then it holds that $z^{\ast}(l) \leq z^{\ast}(f) \leq z^{\ast}(u)$.

\subsection{Bounds for convex cost functions} \label{sec: properties-convex}
\begin{figure}[h]
\begin{subfigure}{0.5\textwidth}
\includegraphics[scale=1]{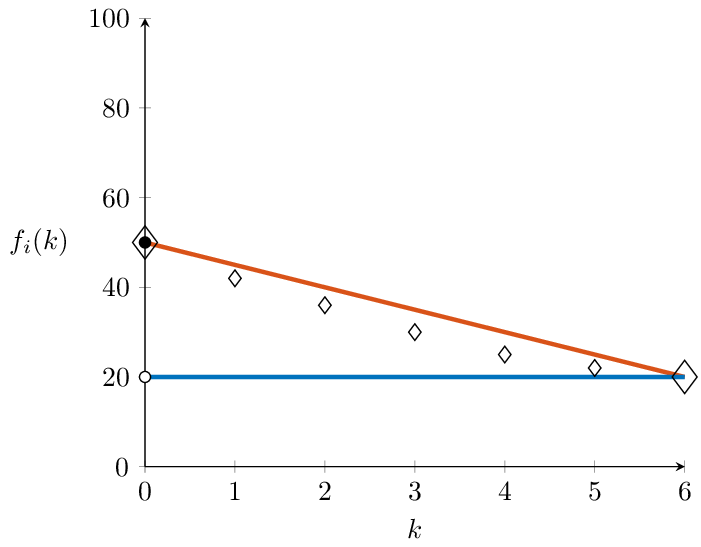}
\caption{\small Two data points \label{fig: piecewise-convex2}}
\end{subfigure}
\begin{subfigure}{0.5\textwidth}
\includegraphics[scale=1]{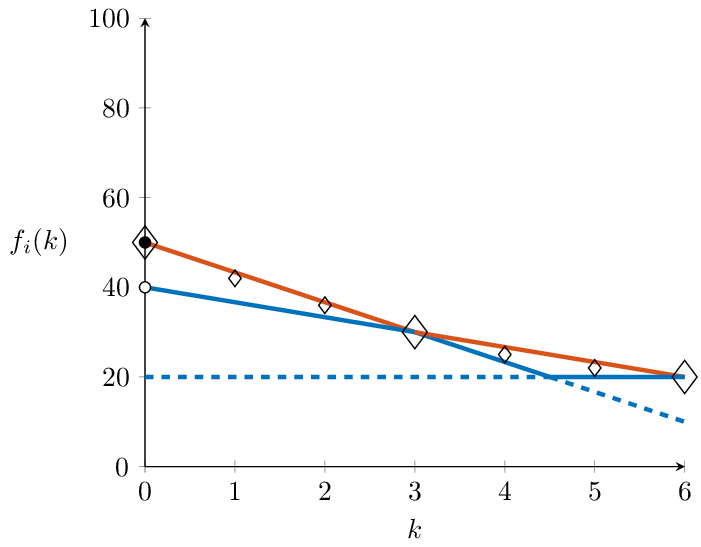}
\caption{\small Three data points. \label{fig: piecewise-convex3}}
\end{subfigure}
\caption{\small Convex cost functions. Evaluated and non-evaluated points are indicated by large and small diamonds, respectively. If $f_i(k)$ is evaluated for only a few values of $k$, we have an upper (red) and lower (blue) bound for the non-evaluated points. The best bounds are indicated by solid lines, redundant bounds are not displayed. Discontinuity points are indicated by circles.\label{fig: piecewise-convex}}
\end{figure}

If cost functions are known to be convex, better lower and upper bounds can be obtained. \Cref{fig: piecewise-convex} gives an example for a single convex (non-increasing) cost function $f_i: \{0,\dotsc,b_i\} \rightarrow [0,M]$, with $b_i=6$. In \Cref{fig: piecewise-convex2} two points are evaluated: $f_i(0)$ and $f_i(6)$. The horizontal blue line illustrates that the value $f_i(6)$ is a lower bound for $f_i(k)$ for all $k$, because $f_i$ is non-increasing. Furthermore, due to convexity of $f_i$ the red line connecting the two points is an upper bound for all $f_i(k)$, $k=0,\dotsc,b_i$. In \Cref{fig: piecewise-convex3} the point $f_i(3)$ is also evaluated. The lines through $f_i(0)$ and $f_i(3)$, and $f_i(3)$ and $f_i(6)$ yield new upper bounds (red) and/or lower bounds (blue) for non-evaluated points $f_i(k)$, $k \in \{0,\dotsc,b_i\}$. The best lower and upper bounds are indicated with solid lines.
We will construct the upper bound value $u_i(k)$ for a particular unobserved point $f_i(k)$; the lower bound is analogous. First, we note that if $v_{i,p}=0$ for all $p=0,\dotsc,b_i$ for the given $i$, then $f_i(k) \leq M$ is an upper bound. Second, if $q>k$ and $v_{i,q} =1$ then
\begin{align}
f_i(k) \leq M \frac{q-k}{q} + f_i(q) \frac{k}{q},
\end{align}
because $f_i(p) \leq M$. Third, if $p<k$ and $v_{i,p}=1$ then  
\begin{align}
f_i(k) \leq f_i(p), 
\end{align}
because $f_i$ is non-increasing. Lastly, if $p<k<q$ and $v_{i,p}=v_{i,q}=1$ then the line connecting these two data points constitutes an upper bound:
\begin{align}
f_i(k) \leq f_i(p)\frac{q-k}{q-p} + f_i(q) \frac{k-p}{q-p},
\end{align}
because $f_i(k)$ is convex. Note that if $p<q<k$ or $k<p<q$ then this line provides a lower bound on $f_i(k)$. Lower bounds are obtained similarly. Similar to \Cref{sec: properties-general}, let $l_i, u_i: \{0,\dotsc,b_i\} \mapsto \mathbb{R}$ denote the functions with the tightest lower and upper bounds for the cost function of player $i$. For evaluated points $f_i(k)$, the lower and upper bounds coincide. Again, let $l(\bm{x}) := \sum_i l_i(x_i)$ and $u(\bm{x}) := \sum_i u_i(x_i)$, and we obtain bounds $z^{\ast}(l) \leq z^{\ast}(f) \leq z^{\ast}(u)$.

\subsection{1-Opt method}
The 1-Opt method can start at any feasible solution, and aims to improve the objective value in each iteration, by moving a single item from one player to another. One beneficial property is that the initial allocation can be varied based on problem specific information.

The 1-Opt method assumes that for each player $i$ two adjacent points on the cost curve are evaluated, and that the set of evaluated points admits an initial allocation $\bm{x}^0$. In each iteration, the method starts with the optimal allocation restricted to only the evaluated data points, and evaluates a single additional point that is adjacent to those already evaluated. The selected new point is the one for which evaluation yields the highest best-case improvement. If the evaluation of such an adjacent point cannot lead to a direct improvement over the current allocation, the method terminates.

Let $\bm{x}^t$ denote the allocation at the start of iteration $t$. For each player $i$, let $\theta_i$ and $\eta_i$ denote the lower bound on the marginal degradation and the upper bound on the marginal gain, respectively:
\begin{subequations}
\begin{align}
%\theta_i^t &= l_{i,x_i^t-1} - f_i(x_i^t) \\
%\eta_i^t &= f_i(x_i^t) - l_{i,x_i^t+1}.
l_i(x_i^t-1) - f_i(x_i^t) \\
\eta_i^t &= f_i(x_i^t) - l_i(x_i^t+1).
\end{align} 
\end{subequations}
In each iteration, $S^{-}$ is the set of players from who an item can be removed, and $\min_{i\in S^{-}} \{\theta_i^t \}$ is the minimum (i.e., best case) cost increase from removing an item from an eligible player. Similarly, $S^{+}$ is the set of players to who an extra item can be allocated, and $\max_{j\in S^{+}} \{\eta_j^t \}$ is the maximum (i.e., best case) cost decrease from adding an item to an eligible player. Let $d^t$ denote the difference:
\begin{align}
d^t := \max_{\substack{ i \in S^{-},j\in S^{+} \\i\neq j}} \{\eta_j^t  - \theta_i^t \}.
\end{align}
As long as $d$ is strictly larger than zero it may be possible to improve the current allocation by moving a single item (i.e., a 1-opt step).

In each iteration $t$, the method considers all players $i \in S^{+}$ with $v_{i,x_i^t +1}^t = 0$, i.e., the players for whom the direct cost decrease after addition of an item is unknown. For such a player $i$, evaluating point $f_i(x_i^t+1)$ yields a best-case improvement of
\begin{align} \label{eq: delta-plus}
\Delta_i^{+,t} = \eta_i^t - \min_j \{\theta_j^t :j \in S^{-}\backslash\{i\},~ v_{j,x_j^t-1}=1\}.
\end{align}
For other players $i$ set $\Delta_i^{+,t} = -\infty$. Allocating an additional item to player $i$ will in the best-case scenario yield an improvement $\eta_i$ for that player. The inner minimization in \eqref{eq: delta-plus} removes the item from the player with the lowest (known) deterioration. Similarly, the method considers all players $i \in S^{-}$ with $v_{i,x_i^t -1}^t = 0$, i.e., the players for who the direct cost increase after removal of an item is unknown. For such a player $i$, evaluating point $f_i(x_i^t-1)$ yields a best-case improvement of
\begin{align} \label{eq: delta-minus}
\Delta_i^{-,t} =\max_j \{\eta_j^t : j \in S^{+}\backslash\{i\} ,~ v_{j,x_j^t+1}=1\} - \theta_i^t,
\end{align}
For other players $i$ set $\Delta_i^{-,t} = -\infty$.

Let $i_t$ be the maximizer\footnote{In case of ties, the player with the lowest index is chosen.\label{tiebreaker-1Opt}} of $\max_i \max\{\Delta_i^{-,t},\Delta_i^{+,t}\}$. Then, either point $f_{i_t}(x_{i_t}^t-1)$ or point $f_{i_t}(x_{i_t}^t+1)$ is evaluated in the current iteration. Subsequently, the optimal allocation on the current set of evaluated points is determined, and the next iteration starts. 

For convex cost functions, the current allocation is changed as follows. If the true cost change (as opposed to the best-case change) is indeed negative, the current allocation is changed according to the 1-opt step. One item is added (if $\Delta_{i_t}^{+,t}>\Delta_{i_t}^{-,t}$) or removed (if $\Delta_{i_t}^{+,t}<\Delta_{i_t}^{-,t}$) from player $i_t$. This item is taken from or moved to player $j_k$, the optimizer\footref{tiebreaker-1Opt} of the inner minimization/maximization of \eqref{eq: delta-plus} or \eqref{eq: delta-minus}. In case of convex cost functions, this 1-opt step is the only change to the allocation in each iteration. On the other hand, for non-convex cost functions it is possible that after a function evaluation there are multiple changes to the current allocation. To ensure the optimal allocation over all currently evaluated points, resource allocation problem \eqref{eq: knapsack} is solved with the additional constraint that only evaluated data points can be used. Pseudocode is given by \Cref{alg: 1-Opt}.

\begin{algorithm}[htb!]
\small
\Begin{
Set $t=0$, set $V^0$ according to initial information, and determine $l^0$\;
Solve $P(f)$ restricted to evaluated data points\;
Denote the solution by $\bm{x}^0$\;
Let $S^{-} = \{i :  x_i^0 >0 \}$ and determine $\theta_i^0$ for all $i\in S^{-}$\;
Let $S^{+} = \{j :  x_j^0 < b_j\}$ and determine $\eta_j^0$ for all $j\in S^{+}$\;
Set $d^0 = \max \{\eta_j^0  - \theta_i^0~|~i \in S^{-},j\in S^{+}, i\neq j\}$\;
\While{$d^t > 0$}{
\For{$i=1:n$}{
	\uIf{$i \in S^{+} \land v_{i,x_i^t+1}^t=0$}{
	$\Delta_i^{+,t} = \eta_i^t - \min_j \{\theta_j^t :j \in S^{-}\backslash\{i\},~ v_{j,x_j^t-1}^t=1\}$\;}
	\Else{$\Delta_i^{+,t} = -\infty$\;}
	\uIf{$i \in S^{-} \land v_{i,x_i^t-1}^t=0$}{
	$\Delta_i^{-,t} = \max_j \{\eta_j^t : j \in S^{+}\backslash\{i\} ,~ v_{j,x_j^t+1}^t=1\} - \theta_i^t$\;}	
	\Else{$\Delta_i^{-,t} = -\infty$\;}
}	
Set $i_t \in \argmax_i \max\{\Delta_i^{-,t}, \Delta_i^{+,t} \}$\;
\uIf{$\Delta_{i_t}^{-,t} < \Delta_{i_t}^{+,t}$}{
Set $\lambda_t = 1$ and let $j_t \in \argmin_j \{\theta_j^t :j \in S^{-}\backslash\{i_t\},~ v_{j,x_j^t-1}^t=1\}$\;
}\Else{
Set $\lambda_t = -1$ and let $j_t \in \argmax_j \{\eta_j^t : j \in S^{+}\backslash\{i_t\} ,~ v_{j,x_j^t+1}^t=1\}$\;
}
Set $V^{t+1} = V^t$\;
Evaluate $f_{i_t}(x_{i_t}^t+\lambda_t)$ and set $v_{i,x_{i_t}^t+\lambda_t}^{t+1}=1$\;
Solve $P(f)$ restricted to evaluated data points\;
Denote the solution by $\bm{x}^{t+1}$\;
Determine $l^{t+1}$\;
Let $S^{-} = \{i :  x_i^{t+1} >0 \}$ and determine $\theta_i$ for all $i\in S^{-}$\;
Let $S^{+} = \{j :  x_j^{t+1} < b_j\}$ and determine $\eta_j$ for all $j\in S^{+}$\;
Set $d^t = \max \{\eta_j^t  - \theta_i^t~|~i \in S^{-},j\in S^{+}, i\neq j\}$\;
Set $t \leftarrow t+1$\;
}
Set $\bm{x}^{\text{1-opt}} = \bm{x}^{t}$\;
}
\caption{1-Opt method \label{alg: 1-Opt}}
\end{algorithm}\vspace*{1em}

The following lemma shows that the 1-Opt method guarantees the optimal solution for convex cost functions.
\begin{lemma} \label{lemma: swapping-convexity}
Let cost functions $f_i: \{0,\dotsc,b_i\} \rightarrow [0,M]$ be convex and non-increasing for all $i=1,\dotsc,n$. Then the 1-Opt solution $\bm{x}^{\text{1-opt}}$ is optimal to \eqref{eq: knapsack}.
\end{lemma}
\begin{proof}
\small
The 1-Opt method terminates with solution $\bm{x}^{\text{1-opt}}$ if and only if $d \leq 0$. This is equivalent to
\begin{align}
% \max_{\substack{ i \in S^{-},j\in S^{+} \\i\neq j}} \Big\{ \big(f_j(x_j^{\text{1-opt}}) - l_{j,x_j^{\text{1-opt}}+1}\big) - \big(l_{i,x_i^{\text{1-opt}}-1} - f_i(x_i^{\text{1-opt}})\big) \Big\} \leq 0.
\max_{\substack{ i \in S^{-},j\in S^{+} \\i\neq j}} \Big\{ \big(f_j(x_j^{\text{1-opt}}) - l_j(x_j^{\text{1-opt}}+1)\big) - \big(l_i(x_i^{\text{1-opt}}-1) - f_i(x_i^{\text{1-opt}})\big) \Big\} \leq 0.
\end{align}
Thus, for each pair $(i,j) \in S^{-}\times S^{+}$ with $i\neq j$ the following inequalities hold:
\begin{align}
\begin{aligned}
f_j(x_j^{\text{1-opt}}) - f_j(x_j^{\text{1-opt}}+1) \leq & 
f_j(x_j^{\text{1-opt}}) - l_j(x_j^{\text{1-opt}}+1) \\
\leq & l_i(x_i^{\text{1-opt}}-1) - f_i(x_i^{\text{1-opt}}) \\
\leq & f_i(x_i^{\text{1-opt}}-1) - f_i(x_i^{\text{1-opt}}),
\end{aligned}
\end{align}
where the first and last inequalities holds because $l_i(k)$ is a lower bound to $f_i(k)$ for all $i,k$. The first term is the `true' cost decrease of assigning an extra item to an eligible player $i$, and the last term is the true cost increase of removing an item from another eligible player $j$. The decrease is smaller than the increase for any player pair $(i,j)$, so moving any item from one player to another cannot yield an improvement to the current allocation. 

Due to convexity of cost functions $f_i$, moving multiple items cannot yield an improvement either. With each additional item involved, the minimum cost increase of removing the item from one player will grow, while the maximum cost decrease of adding the item to another player will diminish. Thus, the current allocation $\bm{x}^{\text{1-opt}}$ is optimal to \eqref{eq: knapsack}. 
\end{proof}
The number of function evaluations depends largely on the set of initially evaluated points. Without problem specific information, one can start with a uniform initialization of matrix $V$, e.g., one can evaluate the points associated with $\floor{B/n}$ and $\floor{B/n}+1$ items for each player (if this is feasible).

The 1-Opt method maintains feasibility in each iteration, and the objective value of the current allocation is always exactly known. Thus, if the method is terminated in any iteration $t$, the current solution $\bm{x}^t$ can directly be implemented. Its objective value is an upper bound for the objective value of the final allocation $\bm{x}^{\text{1-opt}}$. 

\subsection{Sandwich method} \label{sec: sandwich}
The 1-Opt method is a heuristic if the cost functions are not known to be convex. In this section we present a method that yields the optimal solution to the resource allocation problem for both convex and non-convex cost functions. The sandwich method (SW) directly uses the lower and upper bound cost function $l$ and $u$ on the true cost function $f$ that are presented in \Cref{sec: properties-convex} and \Cref{sec: properties-general}. The sandwich method is inspired by \citet{Siem11}, who use sandwich methods, both with and without derivative information, to approximate univariate convex functions.

Let $l^t$ and $u^t$ denote the lower and upper bound cost functions at the start of iteration $t$. In every iteration $t$, we solve the lower bound problem, to obtain objective value $z^{\ast}(l^t)$ and optimal solution $\bm{x}^{l,t}$. The value $z^{\ast}(l^t) = z(\bm{x}^{l,t},l^t)$ is a lower bound to the true optimal objective value $z^{\ast}(f)$. Additionally, $z(\bm{x}^{l,t}, u^t)$ is an upper bound. The sandwich method uses these bounds to formulate an objective value gap. Define
\begin{align}\label{eq: gap}
g^t := z(\bm{x}^{l,t}, u^t) - z(\bm{x}^{l,t}, l^t).
\end{align}
The SW method iteratively improves the lower and upper bounds on the objective value as more data points have been evaluated, and the goal is to evaluate those points that reduce the gap between the upper and lower bound objective value the most. Thus, it `sandwiches' the true cost function (and objective value). As soon as the gap is small enough the method terminates. Solution $\bm{x}^{l,t}$ is the final solution, and the objective value is in the interval $[z(\bm{x}^{l,t}, l^t) , z(\bm{x}^{l,t}, u^t)]$. Note that in any iteration $t$ the value $z(\bm{x}^{u,t}, u^t)$ is also an upper bound on the final objective value, but it is not an upper bound for the objective value corresponding to allocation $\bm{x}^{l,t}$.

As long as the gap $g^t$ is larger than a pre-specified tolerance $\epsilon$, a new iteration starts in which a new point $f_i(k)$ is evaluated; this point is chosen according to some decision rule $\text{DR}$. With $\epsilon=0$ the sandwich method guarantees optimality, both for convex and non-convex cost functions. Pseudocode is given by \Cref{alg: SW}.\\
\begin{algorithm}[H]
\small
\Begin{
Set $t=0$, Pick $\epsilon >0$, set $V^0$ according to initial information\; 
Determine $l^0,u^0$, solve $P(l^0)$ and determine $g^0$\;
\While{$g^t > \epsilon$}{
	Let $(i,k) = \text{DR}(V^t,l^t,u^t)$\;
	Evaluate $f_i(k)$ and determine $V^{t+1}$\;
	Determine cost functions $l^{t+1}$ and $u^{t+1}$\;
	Solve $P(l^{t+1})$ and let $\bm{x}^{l,t+1}$ denote the solution\;
	Determine $g^{t+1}$\;
	Set $t \leftarrow t + 1$\; 	
}
Set $\bm{x}^{\text{sw}} = \bm{x}^{l,t}$\;
}
\caption{Sandwich method \label{alg: SW}}
\end{algorithm}\vspace*{1em}
The corresponding objective value is $z(\bm{x}^{\text{sw}}, f)$. Irrespective of the chosen decision rule, this is within $\epsilon$ of the global optimum.
\begin{lemma} \label{lemma: proof-sandwich}
Let cost functions $f_i : \{0, \dotsc, b_i\} \mapsto [0, M]$ be convex and non-increasing for all $i = 1, \dotsc, n$. Then it holds that $z(\bm{x}^{\text{sw}}, f) \leq z^{\ast}(f) + \epsilon$.
\end{lemma}
\begin{proof}
Let $l^t$ and $u^t$ denote the lower and upper bound matrices at the start of iteration $t$. Let $\bm{x}^{\ast}$, $\bm{x}^{l,t}$ and $\bm{x}^{u,t}$ be minimizers of $P(f)$, $P(l^t)$ and $P(u^t)$, respectively. Then the following sequence of lower bounds holds for $z^{\ast}(f)$:
\begin{align} \label{eq: sandwich-inequalities-LB}
z^{\ast}(l^t) = z(\bm{x}^{l,t}, l^t) \leq z(\bm{x}^{\ast}, l^t) \leq z(\bm{x}^{\ast}, f) = z^{\ast}(f).
\end{align}
The first inequality holds because $\bm{x}^{\ast}$ is feasible for $P(l^t)$ but not necessarily optimal. The second inequality holds because $l_i(k) \leq f_i(k)$ for all $i,k$.  Similarly, the following sequence of upper bounds holds:
\begin{align} \label{eq: sandwich-inequalities-UB}
z(\bm{x}^{l,t}, u^t) \geq z(\bm{x}^{u,t}, u^t)  = z^{\ast}(u^t) \geq z(\bm{x}^{u,t}, f) \geq z(\bm{x}^{\ast}, f) = z^{\ast}(f).
\end{align}
The first inequality holds because $\bm{x}^{l,t}$ is feasible for $P(u^t)$ but not necessarily optimal. The second inequality holds because $u_i(k) \geq f_i(k)$ for all $i,k$, and the third inequality holds because $\bm{x}^{u,t}$ is feasible for $P(f)$ but not necessarily optimal. Thus, at any iteration $t$ the following holds for solution $\bm{x}^{L,t}$:
\begin{align} \label{eq: SW-interval}
z(\bm{x}^{L,t},l^t) \leq z^{\ast}(f) \leq z(\bm{x}^{L,t},u^t).
\end{align}
By construction, if at iteration $t$ the condition
\begin{align} \label{eq: SW-optimality-condition}
z(x_i^{L,t}, f) \leq z^{\ast}(f) + \epsilon,
\end{align}
does not hold, a new function value $f_i(k)$ is evaluated for some $i,k$ in the next iteration in order to reduce the gap. After $n + \bm{b}^{\top}\bm{e}$ iterations, all values $f_i(k)$ have been evaluated for all $k \in \{0,\dotsc,b_i\}$ and all $i=1,\dotsc,n$, so $l_i(k) = f_i(k) = u_i(k)$ for all $i,k$, and the inequalities in \eqref{eq: SW-interval} reduce to equalities. Thus, after at most $n + \bm{b}^{\top}\bm{e}$ iterations the condition \eqref{eq: SW-optimality-condition} must be satisfied.
\end{proof}

It remains to specify the decision rule $\text{DR}$ that picks the new data point to be evaluated in each iteration. Several options are presented below:
\begin{itemize}
\item \textbf{Random (SW-RND)}: Randomly picks a point $(i,k)$ that has not been evaluated yet at the start of iteration $t$, i.e., with $v_{i,k}^t = 0$.

\item  \textbf{Maximum difference - all points (SW-A)}: Determines for all $(i,k)$ the difference between upper bound $u_i^t(k)$ and lower bound $l_i^t(k)$ for the given set of evaluated data points $V^{t}$ at the start of iteration $t$. It then evaluates the data point for which this difference is largest.

\item \textbf{Maximum difference - restricted to points in $\bm{x}^{l,t}$ and/or $\bm{x}^{u,t}$ (SW-R)}: Similar to decision rule SW-A, except that in each iteration only those data points are considered that were chosen in the current allocation $\bm{x}^{l,t}$ and/or $\bm{x}^{u,t}$. If no such points exist, it falls back on SW-A. Note that this rule additionally solves ILP $P(u^t)$ in each iteration.
\end{itemize}
Decision rule SW-R is expected to perform best; SW-RND and SW-A are included in the numerical experiments to illustrate the influence of the employed decision rule.

As a fourth decision rule, we also considered amending SW-R to consider only those data points that were chosen either in the allocation $\bm{x}^{l,t}$ or $\bm{x}^{u,t}$, but not both. The rationale is that for such points, the true value may provide more information than for points that are chosen in both allocations. For the latter, the point is chosen regardless of whether its value is at the lower or upper bound, so allocations are not sensitive to its value. Preliminary experiments indicate that in many iterations $t$ there are no not yet evaluated points that are chosen in $\bm{x}^{l,t}$ or $x^{u,t}$, but not both. Thus, this decision rule often falls back on SW-R, and is not considered further in the numerical experiments in \Cref{sec: numerical-experiments}. 

\begin{figure}[h]
\begin{subfigure}{0.48\textwidth}
\includegraphics[scale=1]{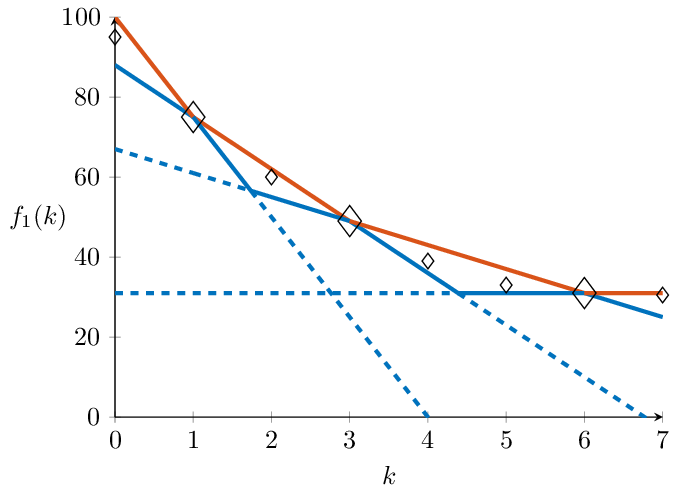}
\caption{\small Player 1. The points $f_1(1)$, $f_1(3)$ and $f_1(6)$ have been evaluated. \label{fig: SW-1}}
\end{subfigure}
\qquad
\begin{subfigure}{0.48\textwidth}
\includegraphics[scale=1]{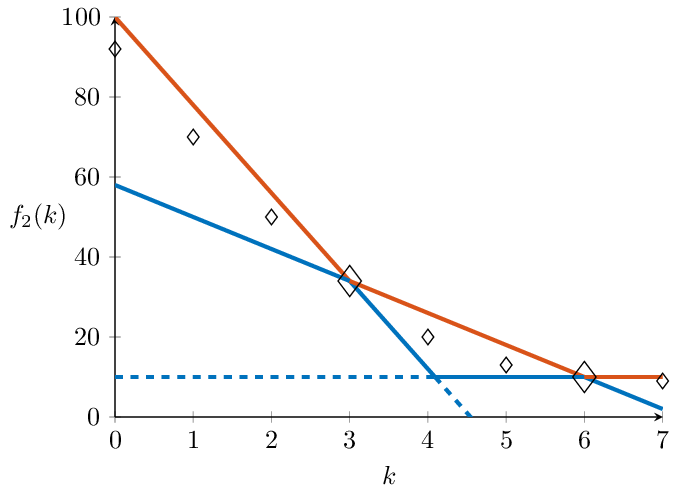}
\caption{\small Player 2. The points $f_2(3)$ and $f_2(6)$ have been evaluated.\label{fig: SW-2}}
\end{subfigure}
\caption{\small Illustration of sandwich method using convexity-based bounds, with allocation budget $B=9$ and upper bound $M=100$. Evaluated and non-evaluated points are indicated by large and small diamonds, respectively. Lower bounds are blue and upper bound are red; the best bounds are indicated by solid lines, redundant bounds are not displayed. The lower bound solution is $\bm{x}^l = (5,4)$ and the upper bound solution is $\bm{x}^u = (4,5)$. \label{fig: SW}}
\end{figure}
\Cref{fig: SW} illustrates the SW method for a situation with two players; for player 1 three points are evaluated and for player 2 two points are evaluated. There are in total 11 non-evaluated points. If decision rule SW-RND is used to determine the next evaluated point, one of these points is selected at random. Decision rule SW-A evaluates point $f_2(0)$, for it has the largest gap. For decision rules SW-R, the solutions to $P(l)$ and $P(u)$ are taken into account. The lower bound solution in the current evaluation is $\bm{x}^l = (5,4)$ and the upper bound solution is $\bm{x}^u = (4,5)$. Among these four points, the gap at $f_2(4)$ is largest: this point is evaluated by SW-R.

\subsection{Illustration of evaluated points}
\begin{figure}[htb!]
\begin{subfigure}{0.48\textwidth}
\hspace*{-0.5cm}
\includegraphics[scale=0.83]{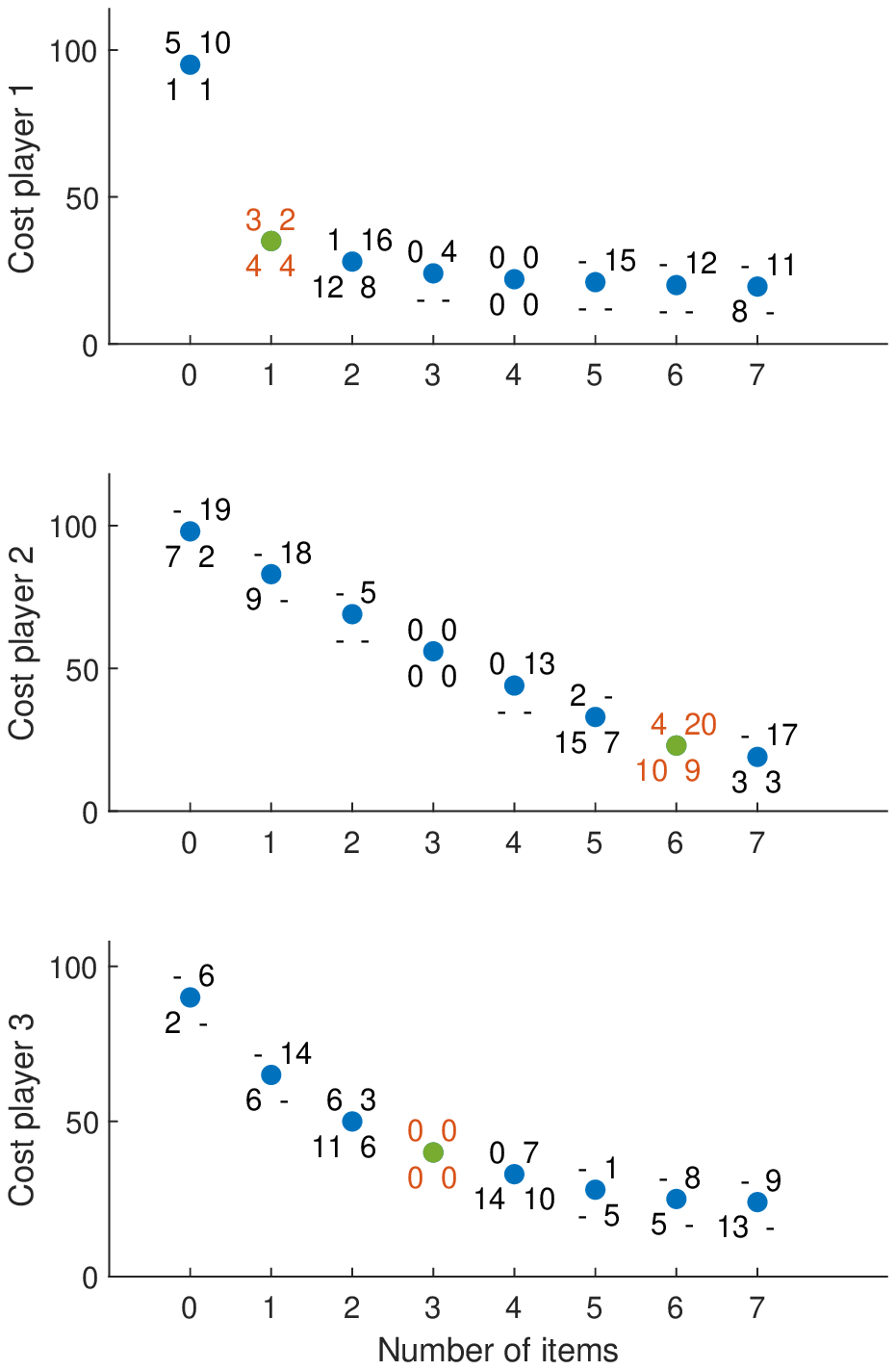}
\caption{\small Convex cost functions.\label{fig: allocation-convex}}
\end{subfigure}
\begin{subfigure}{0.48\textwidth}
\centering
\includegraphics[scale=0.83]{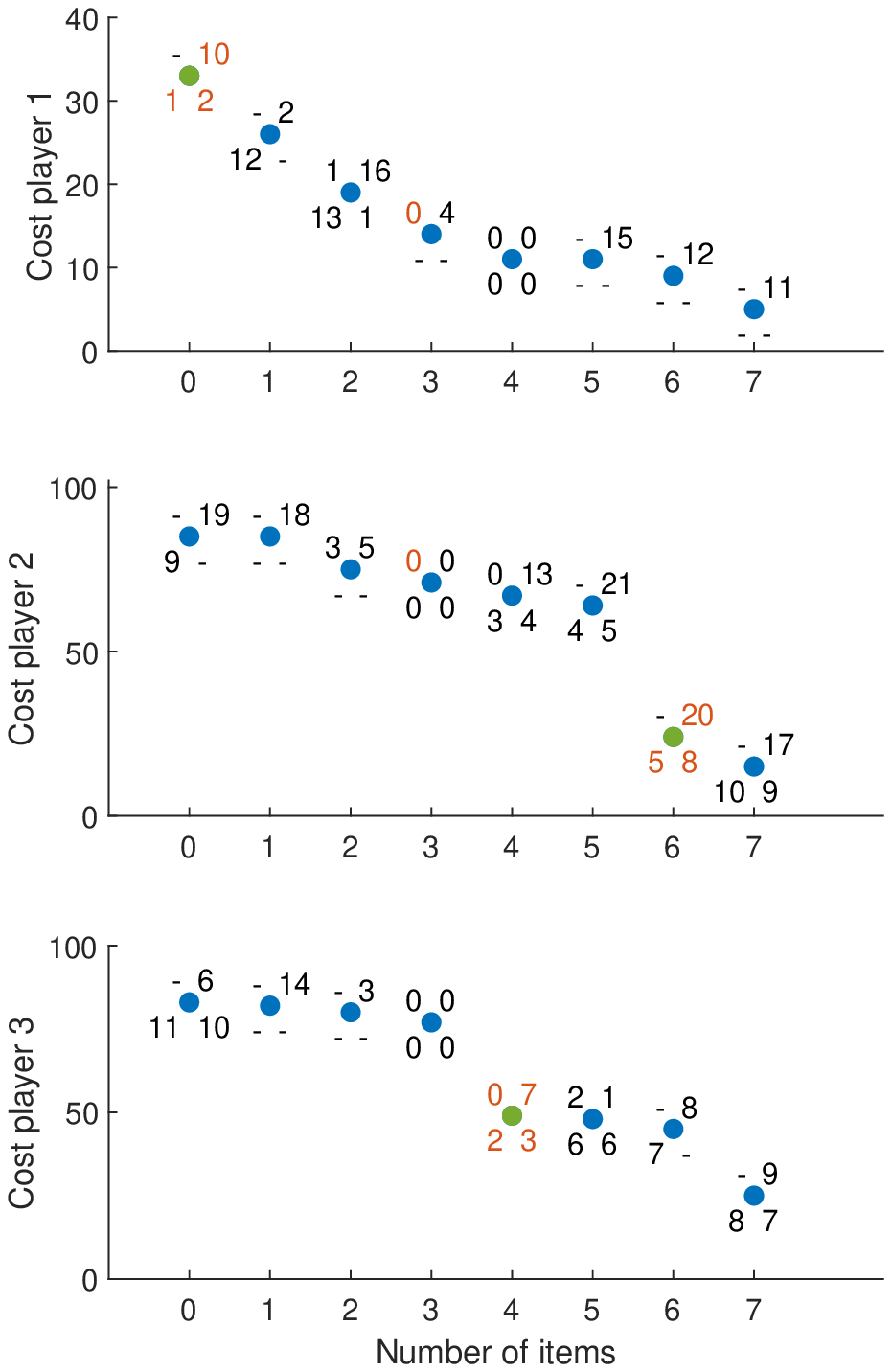}
\caption{\small Non-convex cost functions. \label{fig: allocation-general}}
\end{subfigure}
\caption{Illustration of solution method progress on instances with 3 players, for convex and non-convex cost functions. Per (player, data point), the two numbers indicate in what iteration that point was evaluated for the methods 1-Opt (top left), SW-RND (top right), SW-A (bottom left) and SW-R (bottom right). If the data point is the final allocation for that player for a method, the number is colored red. The symbol `-' indicates that the point was not evaluated by the method. The optimal allocation is indicated with green dots. Vertical axis shifted left for presentation clarity. \label{fig: allocation-example}}
\end{figure}

To compare and contrast the solution methods, \Cref{fig: allocation-example} illustrates the progress of the methods 1-Opt, SW-RND, SW-A and SW-R for two instances with convex and non-convex cost functions ($n=3$, $b_i=7$ for all $i$, $B=10$).
In \Cref{fig: allocation-convex} (convex), the cost function for player 1 has a fast initial drop and marginal improvements after that. For player 2 the improvement rate is near constant, and for player 3 a gradually diminishing improvement can be observed. Overall, 1-Opt, SW-RND, SW-A and SW-R require $12$, $23$, $18$ and $13$ function evaluations, respectively. All methods find the optimal allocation with objective value $98.0$. Method 1-Opt starts with evaluating $f_i(3)$ and $f_i(4)$ for all players. For player 1, it subsequently evaluates data points to the left, whereas for player 2 it evaluates data points to the right of the initial points. For player 3 it evaluates only a single additional point. All sandwich methods start with the initial evaluation of $f_1(4)$, $f_2(3)$ and $f_3(3)$. SW-RND randomly chooses the next points to evaluate, which results in evaluating nearly all points. SW-A and SW-R perform significantly better. SW-A evaluates roughly the same points as 1-Opt. However, it can evaluate non-adjacent points, e.g., for player 2 it evaluates $f_2(k)$ for $k=0,3,5,6,7$.

In \Cref{fig: allocation-general} (non-convex), the cost function for player 1 is non-convex, but decreases gradually (note the difference in vertical axis scaling). The cost functions for players 2 and 3 exhibit drops of larger magnitudes that lead to non-convexity. Method 1-Opt results in an allocation with objective value 134 (9 function evaluations). The sandwich methods all yield the optimal objective value of 106; SW-RND requires 24 function evaluations (i.e., it evaluates all points), SW-A 16 function evaluations and SW-R 13 function evaluations. After a few iterations, 1-Opt cannot find any point where evaluating may result in a direct improvement, and thus terminates. Consequently, it finds only a locally optimal solution.

\subsection{Extensions} \label{sec: extensions}
The problem formulation \eqref{eq: knapsack} can be adapted in various ways. A more general version is obtained by replacing budget constraint \eqref{eq: knapsack-total-budget} by a general linear or convex constraint. Both the 1-Opt and the sandwich method can be applied to such formulations as well, but for ease of exposition we used a simple linear constraint. Another more general version of \eqref{eq: knapsack} is obtained by letting $\varphi : \mathbb{R}^{n} \mapsto \mathbb{R}$ and using the composite objective
\begin{align}\label{eq: composite}
\varphi(f_1(x_1),\dotsc,f_n(x_n)),
\end{align}
where $\varphi(\cdot)$ is an explicitly known function, e.g., the pointwise maximum function. The inequalities \eqref{eq: sandwich-inequalities-LB} and \eqref{eq: sandwich-inequalities-UB} in the proof of \Cref{lemma: proof-sandwich} remain valid, so the sandwich method can be applied to this more general version. The local search concept of the 1-Opt method can be extended to more general composite objective functions as well. However, this requires several modifications, as its current formulation utilizes operations that are specific to the sum function.

\section{Numerical experiments} \label{sec: numerical-experiments}
We test the solution methods on both randomly generated  test instances and instances stemming from applications. \Cref{sec: setup} describes the setup of the numerical experiments: the used benchmark methods and initialization. In \Cref{sec: numerical-results-rnd-convex} and \Cref{sec: numerical-results-rnd-general}, we report and discuss the performance of the solution methods on randomly generated instances with convex and non-convex cost functions, respectively. \Cref{sec: early-termination} considers solution guarantees when the methods are terminated early. In \Cref{sec: RT}, the methods are applied to instances from two applications in radiation therapy treatment planning. In the numerical experiments, problems of form \eqref{eq: knapsack} are solved via their ILP representation (see \Cref{app: ILP}) using Gurobi 9.0 \citep{Gurobi22}.

\subsection{Setup} \label{sec: setup}
To put the performance of the 1-Opt method and the sandwich methods into perspective, we compare the methods with three benchmark methods in terms of their obtained objective value and number of function evaluations. The first benchmark method is the NOMAD solver \citep{NOMAD}, which is the best performing open-source derivative-free solver in our preliminary numerical experiments (see \Cref{sec: literature_review}). NOMAD does not guarantee an optimal solution, neither for the convex nor the non-convex instances. By default, NOMAD does not consider separability of cost functions but only tracks evaluations of the entire cost function $f(\bm{x})$. As such, it does not track function evaluations of cost functions $f_i$ for individual players $i$. In the numerical experiments, this is registered using a custom callback function. However, this does mean that NOMAD may perform iterations without considering new values of $f_i(x_i)$ for any player $i$. Consequently, it does not necessarily recognize that all individual cost functions have been fully evaluated and thus may not terminate at that time.

The other two benchmark methods are simple constructive methods: the myopic method (MY) and the prescient method (PR). The myopic method is a greedy method that starts with zero items allocated, and allocates in each iteration a single item to the player with the largest immediate gain. \Cref{app: myopic} describes the myopic method. The prescient method is similar to the myopic method, but is less greedy. It takes into account the average gain over the remaining horizon for each player, and also assigns items to a player if its immediate gain is low, but its average gain is high. \Cref{app: prescient} describes the prescient method. Both the myopic and the prescient method guarantee the optimal allocation for convex cost functions, but are heuristics for non-convex cost functions.

For the initialization of the 1-Opt method, for each player both the data points corresponding to $\floor{B/n}$ and $\floor{B/n}+1$ items are evaluated to obtain $V^0$. For the sandwich methods we set $\epsilon = 0$, to ensure optimality of the allocations. Matrix $V^0$ is set as follows. We assign $\floor{B/n}$ to each player. The remaining allocation budget $B - n\floor{B/n}$ is evenly allocated, starting from player one; let $\bm{x}^0$ denote the resulting allocation. Matrix $V^0$ is obtained by evaluating all data points corresponding to $\bm{x}^0$.

\subsection{Randomly generated convex cost functions}  \label{sec: numerical-results-rnd-convex}
\begin{figure}[htb!]
\begin{subfigure}{0.46\textwidth}
\includegraphics[scale=1]{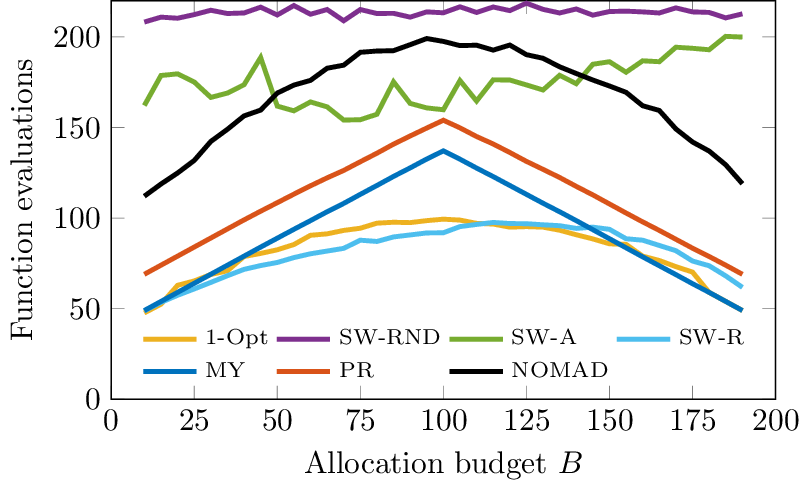}
\caption{\small Small instances. \label{fig: influence_B_convex_small}}
\end{subfigure}
\hspace*{0.04\textwidth}
\begin{subfigure}{0.46\textwidth}
\includegraphics[scale=1]{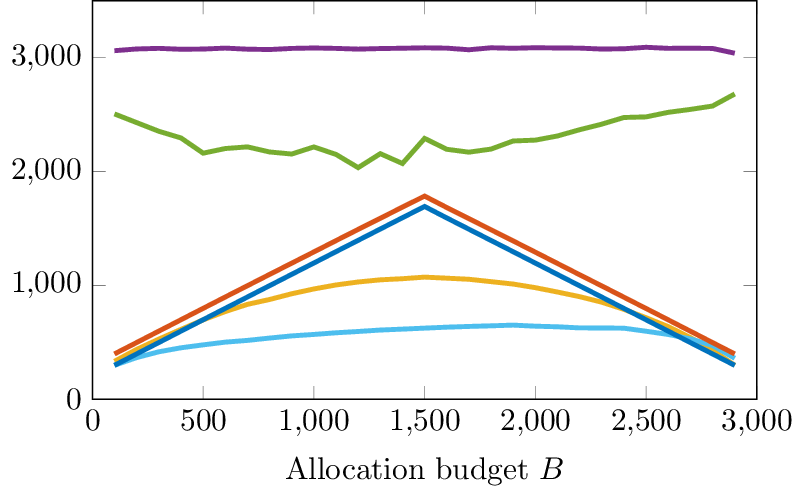}
\caption{\small Large instances. \label{fig: influence_B_convex_large}}
\end{subfigure}
\caption{\small Influence of allocation capacity $B$ on the number of function evaluations. Results of two setups with convex cost functions are shown. The first setup (left panel) uses randomly generated instances with $n=20$ and $b_i=10$ for all $i$, yielding a maximum of 220 function evaluations. The second setup (right panel) uses randomly generated instances with $n=100$ and $b_i=30$ for all $i$, yielding maximum $3{,}100$ function evaluations. Both panels show averages over $K=10$ samples.
\label{fig: influence_B_convex}}
\end{figure}

We compare the performance of the different methods on randomly generated instances with convex cost functions. \Cref{fig: influence_B_convex} compares the required number of function evaluations for the various methods for different choices of allocation budget $B$, for small and large instances. All results are averaged over $K=10$ randomly generated instances.

\Cref{fig: influence_B_convex_small} shows the results for small instances with $n=20$ players and individual budgets of $b_i=10$ for all players, for allocation budget $B$ between $10$ and $190$ (step size 5). The maximum number of function evaluations is $\bm{b}^{\top}\bm{e}+n = 220$.

The sandwich method with decision rules SW-RND and SW-A both perform poorly. Method SW-R requires slightly fewer function evaluations than 1-Opt for low values of $B$, and vice versa for higher values of $B$. We define the average evaluation percentage as follows:
\begin{align}
\text{Evaluation percentage} = \frac{\frac{1}{K}\sum_{k=1}^{K}\# \text{evaluated data points instance } k}{\bm{b}^{\top}\bm{e}+n} \cdot 100\%,
\end{align}
i.e., it is averaged over the $K$ randomly generated instances. Over the entire range of parameter $B$, methods 1-Opt and SW-R have an evaluation percentage of at most $45\%$.

The benchmark methods MY and PR exhibit a very similar pattern. The kink in the curves is where both MY and PR switch from $x^0_i = 0$ for all $i$ to $x^0_i = b_i$ for all $i$ (see \Cref{app: benchmark} for details). PR requires more function evaluations than MY for any $B$, i.e., there is a cost associated to looking beyond the direct gain. The NOMAD solver requires consistently more function evaluations than MY, PR, 1-Opt and SW-R. NOMAD is the only solution method that does not guarantee the optimal solution for convex cost functions; for the experiments of \Cref{fig: influence_B_convex_small} it yields objective values $0.60\%$ higher than optimal, on average.

\Cref{fig: influence_B_convex_large} shows the required number of function evaluations for instances with $n=100$ players and individual budgets of $b_i=30$ for all players. The allocation budget $B$ ranges between $100$ and $2{,}900$ (step size 100). The maximum number of function evaluations is $\bm{b}^{\top}\bm{e}+n = 3{,}100$. The NOMAD solver was not able to obtain a solution for these instances, and is thus excluded from comparison. One possible reason could be that the solution space for NOMAD is vastly larger because it does not exploit separability; whereas NOMAD considers a solution space with $31^{100}$ points, separability-exploiting solvers consider only $3{,}100$ data points, see \Cref{sec: literature_review} for more details.

The results for the small and large instances are comparable for most methods, except for 1-Opt and SW-R. For the larger instances, SW-R requires substantially fewer function evaluations than 1-Opt, unless $B$ is very low or high. One possible reason for this is that with higher individual budgets $\bm{b}$, the distance between the optimal allocation and the initial allocation $\bm{x}^0$ can be larger. Whereas 1-Opt moves only one item per iteration, the sandwich methods can (in theory) change the entire allocation in every iteration. Thus, the sandwich methods may be less influenced by a poor initial allocation. Over the entire range of parameter $B$, methods 1-Opt and SW-R have an evaluation percentage of at most $35\%$ and $21\%$, respectively.

\Cref{table: results-large} shows for each method both the number of function evaluations and the computation time excluding function evaluations, for one large instance with $B=1{,}000$. Under the assumption that each function evaluation is `expensive', the developed methods pose no significant computational burden. 
\begin{table}[h]
\centering
\begin{tabular}{c | c c c c c c c}\toprule
Method & 1-Opt & SW-RND & SW-A & SW-R & MY & PR & NOMAD\\ \midrule
Function evaluations & 1,048 & 3,079 & 2,387 & 565 & 1,193 & 1,285 & - \\
Time (sec) &  19 & 71 & 63 &  25 & 0.6 &  0.9 & -\\ \bottomrule
\end{tabular}
\caption{\small Comparison of all methods on a randomly generated instance with $n=100$, $b_i=30$ for all $i$, $B=1{,}000$ and convex cost functions. Time is the total computation time of the method excluding function evaluations. \label{table: results-large}}
\end{table}

Altogether, the results indicate that, for convex cost functions, the sandwich method with decision rule SW-R requires fewest function evaluations, particularly for large scale instances. For smaller instances, the 1-Opt method is a good alternative. Both SW-R and 1-Opt outperform all benchmark methods.

\subsection{Randomly generated non-convex cost functions}  \label{sec: numerical-results-rnd-general}
\begin{figure}[htb!]
\begin{subfigure}{0.46\textwidth}
\includegraphics[scale=1]{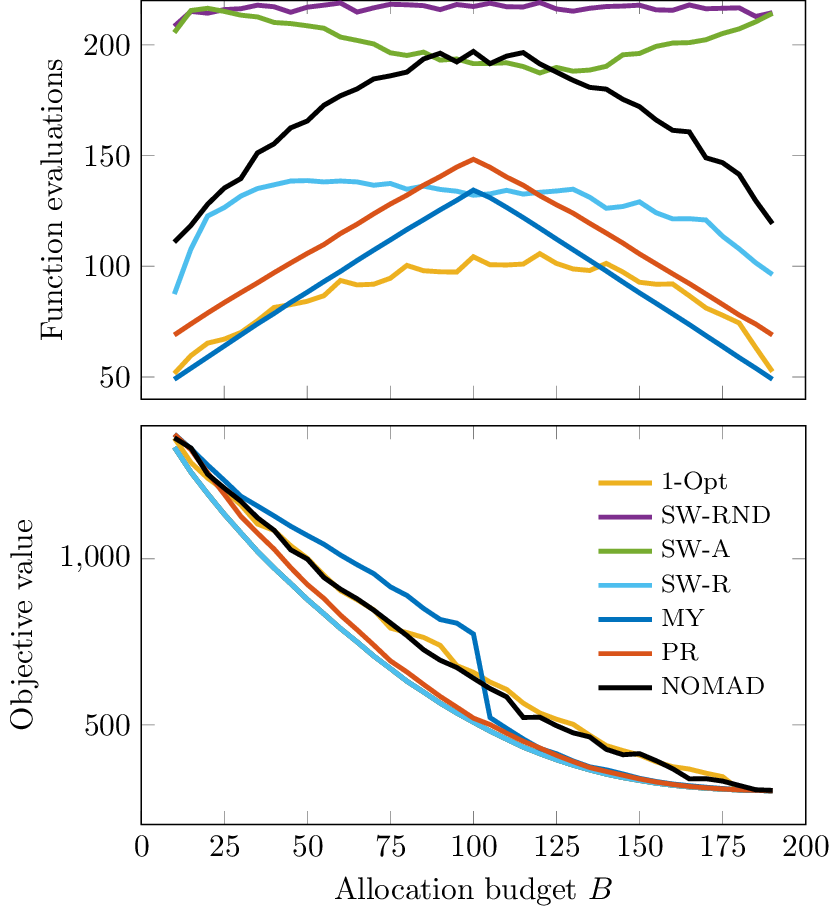}
\caption{\small Small instances. 
\label{fig: influence_B_general_small}
}
\end{subfigure}
\hspace*{0.04\textwidth}
\begin{subfigure}{0.46\textwidth}
\includegraphics[scale=1]{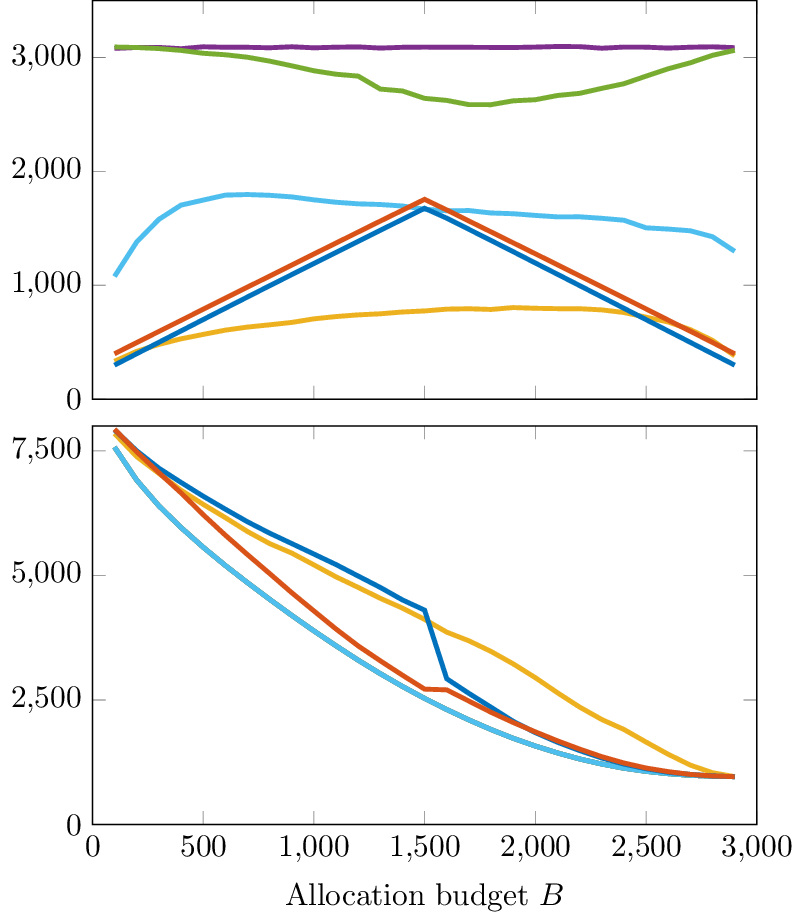}
 \caption{\small Large instances.
 \label{fig: influence_B_general_large}
 }
\end{subfigure}
\caption{\small Influence of allocation capacity $B$ on the number of function evaluations and objective value. Results of two setups with non-convex cost functions are shown. The first setup (left panels) uses randomly generated instances with $n=20$ and $b_i=10$ for all $i$, yielding a maximum of 220 function evaluations. The second setup (right panels) uses randomly generated instances with $n=100$ and $b_i=30$ for all $i$, yielding maximum $3{,}100$ function evaluations. Both panels show averages over $K=10$ samples. In the objective value panels, SW-RND and SW-A coincide with SW-R, which achieves global optimality.
\label{fig: influence_B_general}}
\end{figure}

We compare the performance of the different methods on randomly generated instances with non-convex cost functions. Analogous to \Cref{fig: influence_B_convex}, \Cref{fig: influence_B_general} compares the performance for the various methods for different choices of allocation budget $B$, for small and large instances. All results are again averaged over $K=10$ randomly generated instances.

\Cref{fig: influence_B_general_small} shows the results for small instances with $n=20$ players and individual budgets of $b_i=10$ for all players, for allocation budget $B$ between $10$ and $190$ (step size 5). The maximum number of function evaluations is $\bm{b}^{\top}\bm{e}+n = 220$. Note that, whereas in \Cref{sec: numerical-results-rnd-convex} only NOMAD did not guarantee the optimal solution, for non-convex cost functions 1-Opt and benchmark methods MY, PR are also heuristics. 

In terms of number of function evaluations, the main difference compared to \Cref{fig: influence_B_convex_small} is that SW-R converges slower to the optimal solution. SW-R has an evaluation percentage of at most $63\%$ over the entire range of parameter $B$. Method 1-Opt converges the fastest (i.e., requires fewest function evaluations), but finds poor solutions, on average. Out of the benchmark methods, PR finds solutions with the lowest objective value, at a reasonable number of function evaluations.

\Cref{fig: influence_B_general_large} shows the results for large instances with $n=100$ player, $b_i = 30$ for all players and allocation budget $B$ between $100$ and $2{,}900$ (step size 100). Also here the NOMAD solver was not able to obtain a solution. Again, compared to \Cref{fig: influence_B_convex_large}, SW-R converges slower to the optimal solution. SW-R has an evaluation percentage of at most $58\%$ over the entire range of parameter $B$. Similar to \Cref{fig: influence_B_general_small} the 1-Opt method and the benchmark methods result in suboptimal solutions, with PR being the best performing benchmark method.

The results indicate that for non-convex cost functions, all methods except the sandwich methods are prone to getting stuck in a local minimum. For highly non-convex cost functions, these methods may result in substantially higher objective values.

\subsection{Early termination} \label{sec: early-termination}
A decision maker may decide to terminate the resource allocation procedure early, perhaps due to financial or time restrictions, with the aim to obtain an implementable allocation. For the 1-Opt method, the algorithm can be terminated in any iteration to obtain a feasible solution. Additionally, by construction, all objective values corresponding to the allocation have been computed, so the objective value is known exactly. For the sandwich method, termination in any iteration also results in a feasible solution, but the objective value is not necessarily precisely known. However, the objective value does satisfy inequalities \eqref{eq: SW-interval}. In this section, we report results on the objective value convergence of 1-Opt method and the sandwich method.

\begin{figure}[htb!]
\begin{subfigure}{0.48\textwidth}
\includegraphics[scale=1]{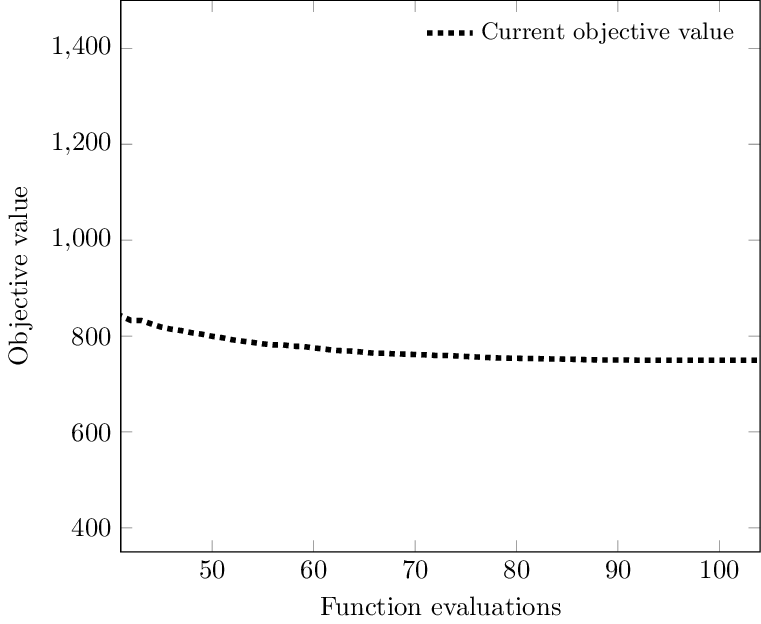}
\caption{\small 1-Opt \label{fig: convergence-convex-1-Opt}}
\end{subfigure}
\hspace*{0.02\textwidth}
\begin{subfigure}{0.48\textwidth}
\includegraphics[scale=1]{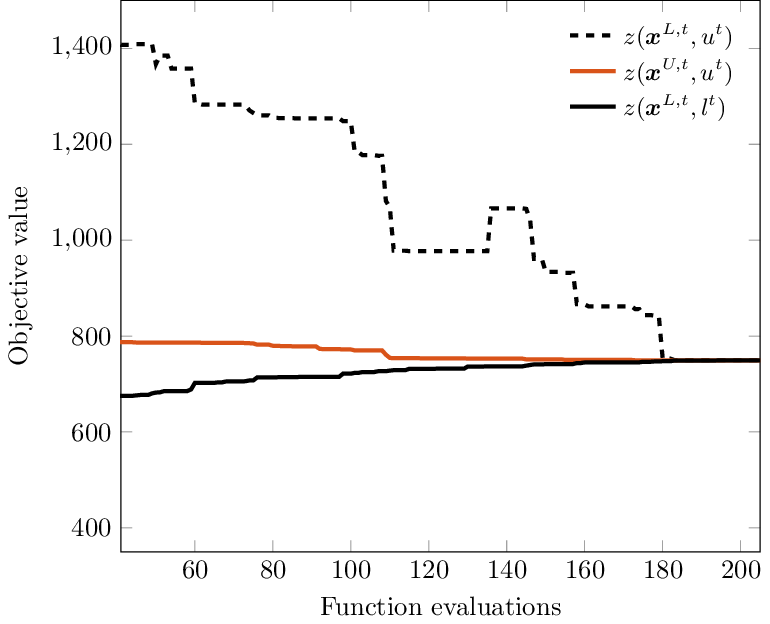}
\caption{\small SW-RND \label{fig: convergence-convex-SW-RND}}
\end{subfigure}
\begin{subfigure}{0.48\textwidth}
\includegraphics[scale=1]{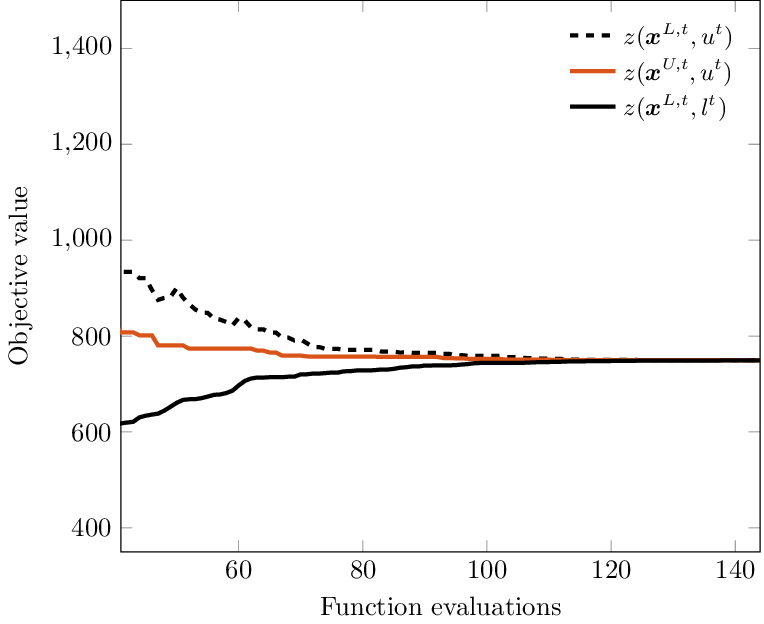}
\caption{\small SW-A \label{fig: convergence-convex-SW-A}}
\end{subfigure}
\hspace*{0.02\textwidth}
\begin{subfigure}{0.48\textwidth}
\includegraphics[scale=1]{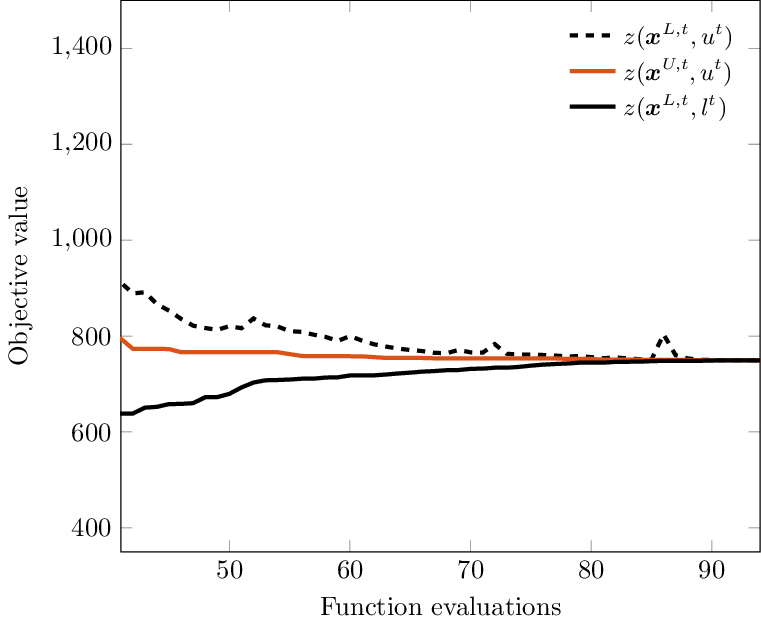}
\caption{\small SW-R \label{fig: convergence-convex-SW-R}}
\end{subfigure}
\caption{\small Convergence of the 1-Opt and sandwich methods on a randomly generated instance with $n=20$, $b_i=10$ for all $i$, $B=90$ and convex cost functions. Note the different scaling on the horizontal axis. Function evaluations during initialization are not displayed, and each method terminates at the maximum displayed value on the horizontal axis. \label{fig: convergence-convex}}
\end{figure}

For the sandwich method, the solution $\bm{x}^{L,t}$ yields at least objective value $z(\bm{x}^{L,t},l^t)$ and at most objective value $z(\bm{x}^{L,t},u^t)$. The upper bound $z(\bm{x}^{L,t},u^t)$ is not necessarily monotonically decreasing in each iteration, because $\bm{x}^{L,t}$ is feasible but not necessarily optimal to $P(u^t)$. In each iteration $t$ one can additionally solve problem $P(u^t)$ to obtain solution $\bm{x}^{U,t}$, with objective value $z(\bm{x}^{U,t},u^t)$. This is also an upper bound to the optimal objective value $z^{\ast}(f)$, although not necessarily an upper bound to the objective value associated with solution $\bm{x}^{L,t}$. In the subsequent result, we report both upper bounds.

\Cref{fig: convergence-convex} shows the convergence of the 1-Opt and sandwich methods on one randomly generated instance with convex cost functions, with $n=20$, $b_i=10$ for all $i$ and $B=90$. The function evaluations during the initialization step are excluded, because they are not informative. For the sandwich methods, a substantial number of function evaluations is required to reduce the gap(s) to exactly zero. For 1-Opt, the objective value changes only slightly in the last 20 iterations. This suggests that for applications where good but not necessarily optimal allocations are required, early termination can be beneficial.

\begin{figure}[htb!]
\begin{subfigure}{0.48\textwidth}
\includegraphics[scale=1]{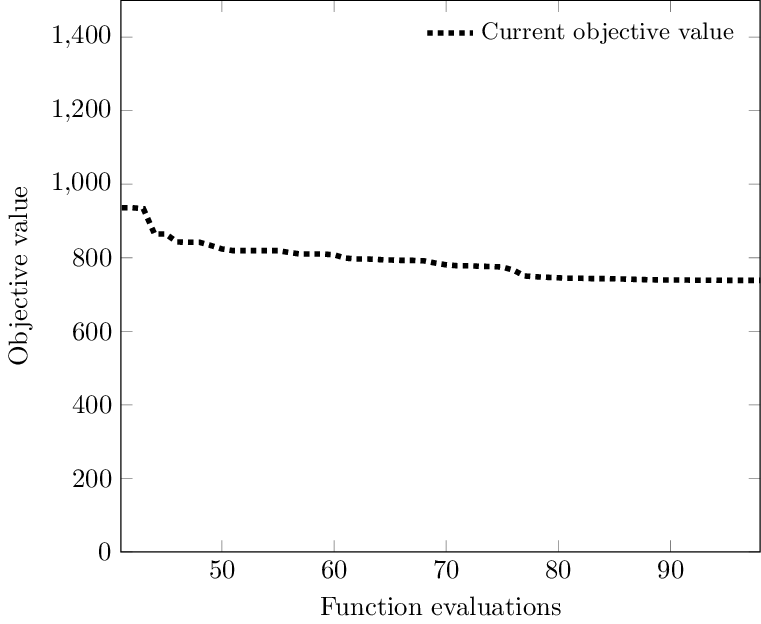}
\caption{\small 1-Opt ($738$) \label{fig: convergence-general-1-Opt}}
\end{subfigure}
\hspace*{0.02\textwidth}
\begin{subfigure}{0.48\textwidth}
\includegraphics[scale=1]{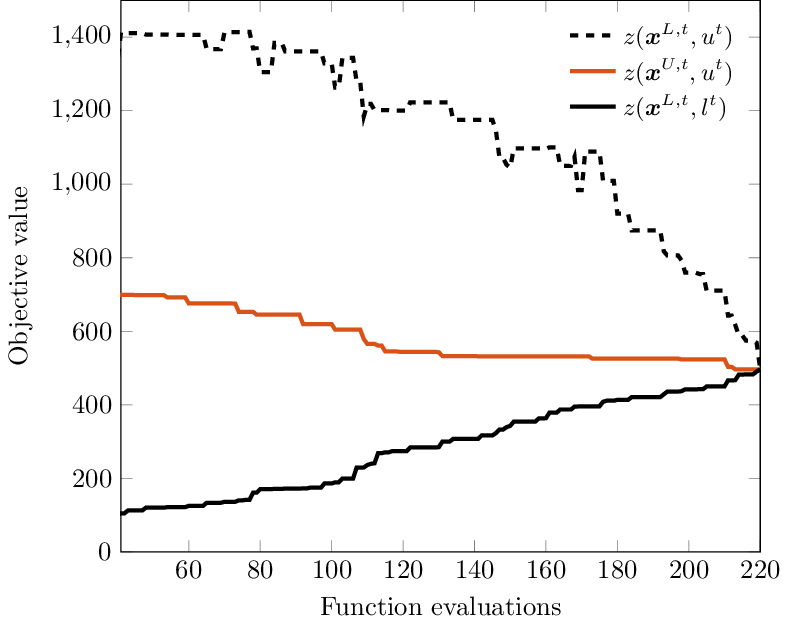}
\caption{\small SW-RND ($497$)\label{fig: convergence-general-SW-RND}}
\end{subfigure}
\begin{subfigure}{0.48\textwidth}
\includegraphics[scale=1]{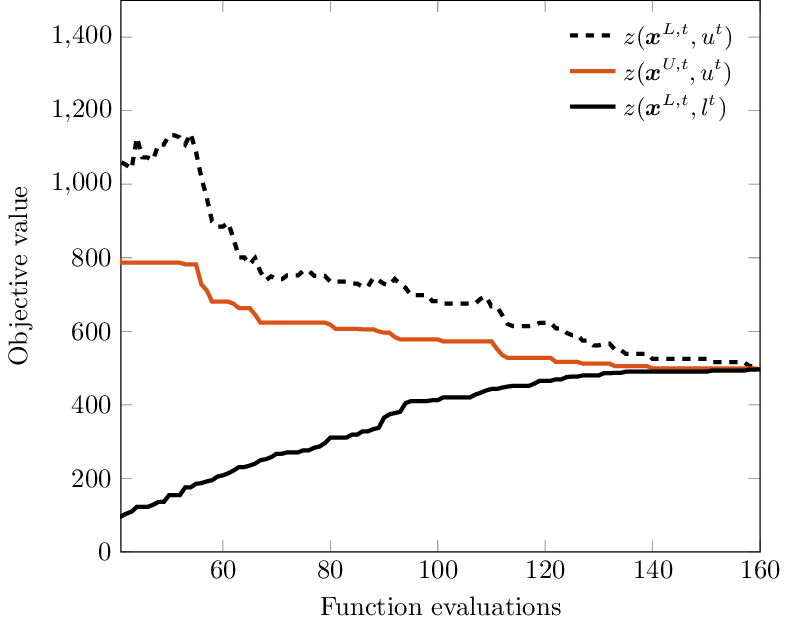}
\caption{\small SW-A ($497$) \label{fig: convergence-general-SW-A}}
\end{subfigure}
\hspace*{0.02\textwidth}
\begin{subfigure}{0.48\textwidth}
\includegraphics[scale=1]{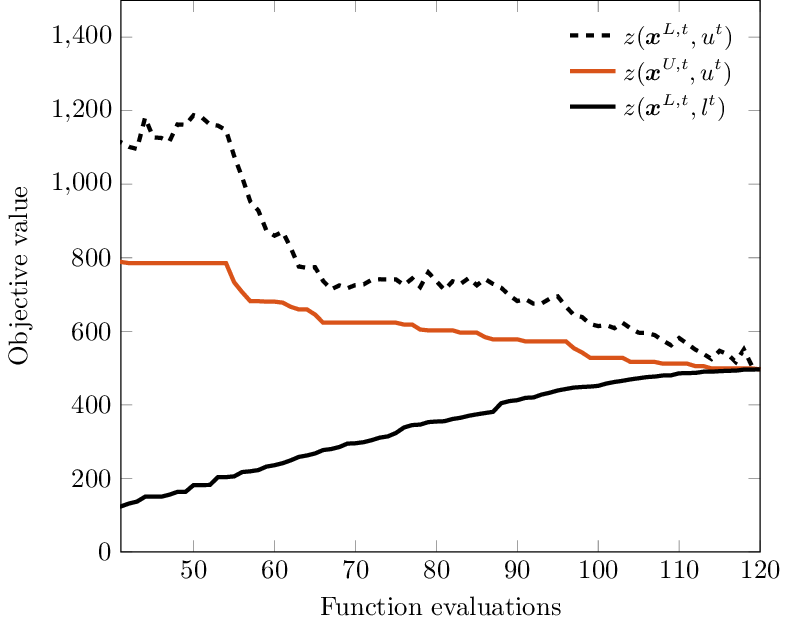}
\caption{\small SW-R ($497$)  \label{fig: convergence-general-SW-R}}
\end{subfigure}
\caption{\small Convergence of the 1-Opt and sandwich methods on a randomly generated instance with $n=20$, $b_i=10$ for all $i$, $B=90$ and non-convex cost functions. Note the different scaling on the horizontal axis. Function evaluations during initialization are not displayed, and each method terminates at the maximum displayed value on the horizontal axis. The resulting objective value is indicated in brackets. \label{fig: convergence-general}}
\end{figure}

\Cref{fig: convergence-general} shows the convergence of the 1-Opt and sandwich methods on one randomly generated instance with non-convex cost functions, with $n=20$, $b_i=10$ for all $i$ and $B=90$. Similar to \Cref{fig: convergence-convex}, the function evaluations during the initialization step are excluded. The main observation is that, compared to \Cref{fig: convergence-convex}, for the sandwich methods the gaps between the lower and upper bound(s) are larger. There is also a considerable difference between both upper bounds. The latter is indicative of the uncertainty of the solution in the current iteration: upper bound $z(\bm{x}^{L,t}, \bm{u}^t)$ is the upper bound for the current solution, but the upper bound for the optimal solution, $z(\bm{x}^{U,t}, \bm{u}^t)$, is much lower. Even in the last few iterations the gaps can be considerable, so early termination of the methods can result in large objective value uncertainty. This example also illustrates that for non-convex cost functions 1-Opt can get stuck in a local optimum; it performs the lowest number of function evaluations (100), but the objective value (738) is substantially higher than the optimal objective value (497).

\subsection{Two applications in radiation therapy planning} \label{sec: RT}
Radiation therapy (RT) is one of the predominant treatment modalities for cancer, and treatment planning for RT relies heavily on mathematical optimization. In the current section, we describe two applications in RT treatment planning that give rise to resource allocation problems with expensive function evaluations, and illustrate the developed methods on instances of both applications. 

In both applications, the allocation cost functions need not be convex in theory, but they are near-convex in practice. This means that 1-Opt and the benchmark methods MY, PR do not guarantee optimality. For 1-Opt and the sandwich methods, we report results using bounds with and without the convexity assumption. Note that if the convexity assumption is used while cost functions are in fact not convex, the sandwich methods lose their optimality guarantee. However, with the convexity assumption tighter lower and upper bounds are obtained, which might translate to fewer function evaluations.

\subsubsection*{Proton therapy slot allocation} 
Proton therapy is an advanced RT modality, which, in theory, offers several advantages over conventional photon therapy (i.e., X-rays). For many patients, proton therapy treatments can deliver more radiation to the tumor, with better sparing of surrounding healthy tissues. As a result, patients may exhibit fewer negative radiation-induced side-effects. There is a high demand for proton therapy, but its availability is scarce due to its high cost. There is a limited number of hospitals that offer proton therapy, and at these facilities the capacity is limited. Thus, hospitals that offer proton therapy try to use their available capacity as efficient as possible. For some patients, the benefits of proton therapy over photon therapy are more pronounced than for others. In the Netherlands, a model-based approach is used to determine who is eligible for proton therapy \citep{Langendijk13}. Recently, research also focusses on delivering only part of the treatment with proton therapy \citep{tenEikelder19a}. We investigate the proton therapy resource allocation problem for such combined proton-photon treatments, which has previously been discussed by \citet{Loizeau21} for a situation without expensive function evaluations.

The treatment for each patient is split up in multiple sessions, known as fractions. Part of the treatment may be delivered using photon therapy; such treatment slots are assumed to be widely available. The algorithm developed in \citet{tenEikelder19a} determines the optimal combined proton-photon treatment schedule for a patient, i.e., the optimal combination of proton and photon fractions (and their treatment dosage), given constraints on these resources. The corresponding objective value is the mean biological effective dose \citep[BED,][]{Fowler89,Hall12} to the tumor that can be attained using this treatment schedule (yielding tolerable doses to surrounding healthy tissues). Let $\text{BED}_i(p)$ denote the mean tumor BED that can be attained if the treatment for patient $i$ can use at most $p$ proton slots. Varying $p$ allows us to quantify the gain of allocating patient $i$ a particular number of proton therapy slots, but each evaluation of $\text{BED}_i(p)$ requires a computationally demanding run of the algorithm in \citet{tenEikelder19a}.

Consider a proton treatment facility with a limited number of proton slots $B$ that are available during a given time horizon, and a set of $n$ cancer patients that are scheduled for treatment in this time horizon and might qualify for proton therapy as (part of) their treatment. Each patient $i$ may be treated with at most $b_i$ proton fractions. The goal is to allocate the proton slots over the set of patients such that the sum of individual BED values is maximized. Let $x_i$ be decision variables, i.e., the number of proton slots allocated to patient $i$. The optimization model reads
\begin{subequations}\label{eq: proton-slot}
\begin{align}
\max_{\bm{x}}~&~ \sum_{i=1}^n \text{BED}_i(x_i) \label{eq: proton-slot-1}\\
\text{s.t.} ~&~ \sum_{i=1}^n x_i = B \label{eq: proton-slot-2}\\
~&~ 0 \leq x_i \leq b_i,~x_i \text{ integer},~\forall i=1,\dotsc,n.
\end{align}
\end{subequations}
We consider an instance with $n=17$ and $b_i=15$ for all $i$, based on real patient data of liver cancer patients treated with proton therapy at Massachusetts General Hospital (Boston, USA). The BED functions for all patients are plotted in \Cref{fig: proton-slot-payoff} in \Cref{app: RT-data}; note that they are near-concave. Evaluation of a single point on a curve takes on average 3 minutes. To transform \eqref{eq: proton-slot} to a minimization problem of form \eqref{eq: knapsack}, function $f_i: \{0,\dotsc,b_i\} \mapsto \mathbb{R}$ is defined as $f_i(k) = 200 - \text{BED}_i(k)$ for each $i,k$.\footnote{For ease of exposition, we use range bounds $[0,M]$. For this problem instance, 200 is an upper bound to $BED_i(k)$ for each player $i$; we add $M=200$ so that every function has a range of $[0,M]$.} We assume that the total number of available proton therapy slots is $B=120$; the sum of individual budgets $b_i$ equals $255$.

The sandwich methods and 1-Opt are initialized according to \Cref{sec: setup}, for MY and PR see \Cref{app: benchmark}. \Cref{table: proton-results} shows the results (of the minimization problem). All methods find the optimal solution (objective value 1,666) or a near-optimal solution. Method 1-Opt requires a low number of function evaluations, both with convex and non-convex lower and upper bounds. It yields the optimal objective value in both cases, although there is no theoretical guarantee for this. Unlike in \Cref{fig: convergence-general} it does not get stuck in a local minimum; the BED functions shown in \Cref{fig: proton-slot-payoff} are `almost' concave (i.e., the cost functions are almost convex). With non-convex cost functions, the sandwich methods are guaranteed to find the optimum, but require many iterations. Particularly for SW-R, the number of iterations decrease considerably when bounds based on convexity are used, with only a small deterioration in objective value. Decision rules SW-RND and SW-A perform poorly in both cases.

Thus, SW-R with convexity based bounds and 1-Opt are the best performing solution methods. These also outperform the benchmark methods: MY yields a near-optimal solution in 151 function evaluations. Method PR uses more function evaluations than MY, but this does not result in a better objective value. NOMAD requires a large number of function evaluations.
\begin{table}[h]
\centering
\begin{tabular}{c c | c c }\toprule
Method& Bound type & Function evaluations & Objective value \\\midrule
\multirow{2}{*}{1-Opt} & Convex & 109 & 1,666\\
& Non-convex & 117 & 1,666\\[0.5em]
\multirow{2}{*}{SW-RND} & Convex & 242 & 1,666\\
& Non-convex & 257 & 1,666\\[0.5em]
\multirow{2}{*}{SW-A} & Convex & 198 & 1,668\\
& Non-convex & 262 & 1,666\\[0.5em]
\multirow{2}{*}{SW-R} & Convex & 109 & 1,668 \\
& Non-convex & 242 & 1,666\\[0.5em]
MY & - 						& 151 & 1,672\\
 PR & -						& 166 & 1,672 \\
NOMAD	&-	& 222 & 1,666 \\ \bottomrule
\end{tabular}
\caption{\small Results for the proton slot allocation problem (transformed to a minimization problem). Bound type refers to the use of lower and upper bounds with or without the convexity assumption. The maximum number of function evaluations is 272. Computation time (excluding function evaluations) for NOMAD was 20 minutes; for all other methods less than 5 seconds. \label{table: proton-results}}
\end{table}

As an alternative objective, one can maximize, e.g., the minimum BED over the patient cohort. This can be modelled by using a different function $\varphi(\cdot)$ in the composite objective \eqref{eq: composite} in \Cref{sec: extensions}.

\subsubsection*{Volumetric modulated arc therapy} 
Whereas the first application considers allocation of resources to patients, the second application consider an allocation problem within the treatment planning process. The resource that has to be allocated efficiently is time itself.

During RT delivery, the patient is positioned on a treatment couch, while a treatment gantry rotates around the patient to deliver radiation from different angles. Using different angles leads to better sparing of healthy tissues surrounding the tumor. Volumetric modulated arc therapy (VMAT) is an advanced RT technique in which radiation is delivered continuously, while the treatment gantry rotates around the patient \citep{Otto08}. 

Planning methods typically split up the full gantry arc into a finite set of arc segments \citep{Balvert17b,vanAmerongen17}. From the beam's eye view, the geometry of tumor and surrounding healthy tissues is two-dimensional, and is different per arc segment. Thus, the desired radiation `pattern', referred to as the fluence map, differs per arc segment. By dynamically blocking parts of the beam during irradiation using a multi-leaf collimator (MLC), the desired fluence map can be approximated. For further details on this procedure and other parts of treatment plan optimization we refer to \citet{Ehrgott08}. 

The quality of approximation of the desired fluence map in an arc segment is higher if more time is allocated to treatment delivery in that arc segment. If sufficient time is allocated to all treatment arcs, all fluence maps can be perfectly replicated. However, such long treatments are undesirable as the expected inaccuracy caused by patient movement increases, and put a strain on hospital facilities \citep{Kelly19}. It is common to aim for treatments with a low delivery time, e.g., at most $T_{\text{total}}$ per patient. For VMAT treatments, this time has to be allocated efficiently over the set of arc segments. 

Time is discretized into time steps of $\Delta$ seconds. Let $M$ denote the number of fluence maps (one per arc segment). Replicating a single fluence map $m=1,\dotsc,M$ for a given arc segment and a particular number of allocated time steps is a non-convex optimization problem in itself and is computationally expensive \citep{vanAmerongen17}. Let $h_m: \mathbb{N} \mapsto \mathbb{R}+$ denote the fluence map matching inaccuracy (in sum-of-squared differences (ssdif)) for a given amount of time steps. Any arc segment $m$ must be allocated between $1$ and $T_m$ time steps. Let $x_m$ be the decision variable, i.e., the number of time steps allocated to arc segment $m$. The allocation problem reads
\begin{subequations}
\begin{align}
\min_{\bm{x}} ~&~ \sum_{m=1}^M h_m(\Delta x_m) \\
\text{s.t.} ~&~ \sum_{m=1}^M x_m \leq T_{\text{total}}/\Delta \\
~&~ 1\leq x_m \leq T_m,~x_m \text{ integer},~m=1,\dotsc,M.
\end{align}
\end{subequations}
Evaluation of $h_m(\Delta t)$ for some $t$ is time expensive, and, in order to facilitate fast treatment planning, the number of function evaluations should be minimized. Note that in practice, evaluation time of $h_m(\Delta t)$ increases for larger $t$, so minimizing the number of function evaluations is an approximation of the `true' objective of minimizing function evaluation time.

We consider the prostate cancer instance of the open-source CORT data set \citep{Craft14}. The full treatment arc is split into $M=36$ arc segments with associated fluence maps. For each arc segment $m$, the maximum time $T_m$ ranges between $5$ and $8$ seconds, allocating more time to the arc segment will not improve the fluence map matching. The time step size is $\Delta = 1$ second. The fluence map matching inaccuracy functions for all arc segments are plotted in \Cref{fig: VMAT_data} in \Cref{app: RT-data}; note that they are near-convex. Evaluating a single point on a curve takes on average 2 minutes. We assume that the entire treatment cannot take more than $T_{\text{total}} = 165$ seconds; the sum of individual bounds $T_m$ equals $254$.

\begin{table}[h] %Original (uniform) starting point y0
\centering
\begin{tabular}{c l l | c c }\toprule
Method & Bound type & Starting point & Function evaluations & Objective value \\\midrule
\multirow{4}{*}{1-Opt} & Convex &Uniform & 112 & 1,072\\%For starting point: B/n and y0
& Non-convex& Uniform & 113 & 970\\ %For starting point: B/n and y0
& Convex &Informed & 99 & 1,072\\ %Informed starting point based on V0
& Non-convex &Informed & 97 & 970\\[0.5em] %Informed starting point based on V0
\multirow{4}{*}{SW-RND} & Convex &Uniform & 247 & 1,013\\%For starting point: B/n and y0
& Non-convex &Uniform & 254 & 950\\ %For starting point: B/n and y0
& Convex &Informed & 251 & 970\\%Informed starting point based on V0
& Non-convex &Informed & 252 & 950\\[0.5em]%Informed starting point based on V0
\multirow{4}{*}{SW-A} & Convex &Uniform & 182 & 1,188\\ %For starting point: B/n and y0
& Non-convex& Uniform & 186 & 950\\ %For starting point: B/n and y0
& Convex &Informed & 193 & 1,028\\%Informed starting point based on V0
& Non-convex &Informed & 184 & 950\\[0.5em]%Informed starting point based on V0
\multirow{4}{*}{SW-R} & Convex &Uniform & 140 & 1,033 \\ %For starting point: B/n and y0
& Non-convex& Uniform & 167 & 950\\ %For starting point: B/n and y0
& Convex &Informed & 115 & 1,072 \\ %Informed starting point based on V0
& Non-convex& Informed & 160 & 950\\[0.5em]%Informed starting point based on V0
MY & -		&-		& 160 & 970\\
 PR & -		&-		& 193 & 970\\
\multirow{2}{*}{NOMAD} &  -& Uniform & 246 & 1,066\\
  & - & Informed &  236 & 984\\ \bottomrule
\end{tabular}
\caption{\small Results for the VMAT problem, both with a uniform starting point and an informed starting point. The maximum number of function evaluations is 254. Computation time (excluding function evaluations) for all methods except NOMAD is less than 5 seconds; NOMAD was terminated manually after 1 hour, for the reason described in \Cref{sec: setup}. \label{table: VMAT-results}}
\end{table}
The performance of 1-Opt and the sandwich methods depends on the starting point. We compare two starting points. The initialization of \Cref{sec: setup} is referred to as the `uniform starting point'. As an alternative, we use application specific information to construct a more informed starting point. In VMAT, it is possible to use simple heuristics to estimate the `toughness' of replicating a certain fluence map \citep{vanAmerongen17}. Matrix $V^0$ is obtained by evaluating two points per curve; for fluence maps with a high toughness measure the evaluated points correspond to more time steps. This is referred to as the `informed starting point', and is used for both 1-Opt and the sandwich methods.

\Cref{table: VMAT-results} shows the results, for both starting points and for both bounds with and without the convexity assumption. 1-Opt does not find the optimum, but is near-optimal for the cases using bounds without the convexity assumption. For 1-Opt, the informed starting point results in a lower number of function evaluations than the uniform starting point, without a change in the resulting objective value, both with and without the convexity assumption. 1-Opt requires substantially fewer function evaluations than the other methods in all situations. This is in contrast with the proton slot allocation problem. There, using the convexity assumption, SW-R required the same number of function evaluations as 1-Opt. For the current application, the convexity assumption results in substantially fewer function evaluations for SW-R, but still more than 1-Opt. When using bounds based on the convexity assumption, the sandwich methods yield suboptimal objective values for both starting points. Without the convexity assumption, the sandwich methods yield the optimal objective value by construction. 	

Both benchmark methods MY and PR result in the same objective value as 1-Opt, while requiring more function evaluations. Similar to the proton slot allocation problem, PR uses more function evaluations than MY, but this does not translate to a better objective value. NOMAD again requires more function evaluations than MY and PR, while resulting in a worse objective value. SW-R without the convexity assumption requires a similar number of function evaluations as MY, but it yields a better objective value.

The VMAT results show that both the choice between bounds with or without the convexity assumption and the choice of starting point can influence the resulting objective value and number of function evaluations. 1-Opt can yield near-optimal solutions while requiring substantially fewer function evaluations than the sandwich methods.

\section{Conclusion} \label{sec: conclusion}
In practical applications of resource allocation problems, evaluating the objective function for a given number of resources can be expensive. Thus, it may be necessary to use only few function evaluations during the decision-making process. However, current solution methods for resource allocation problems do not take this aspect into account. At the same time, derivative-free integer optimization methods do not properly use the separable objective structure that is typically encountered in resource allocation problems. To bridge this gap, we introduced solution methods for integer resource allocation problems that aim to limit the number of function evaluations. We have considered formulations with both convex and non-convex separable cost functions; for convex cost functions both the 1-Opt and sandwich methods guarantee an optimal solution. For non-convex cost functions, only the sandwich methods guarantee optimality. 

From our numerical experiments, we conclude that for convex cost functions the sandwich method SW-R requires fewest function evaluations on large-scale instances, and 1-Opt is comparable to SW-R for small instances. For non-convex cost functions, 1-Opt uses fewer function evaluations than SW-R, but is prone to getting stuck in a local minimum. This may result in substantially worse allocations than the sandwich methods (which are exact by construction). Depending on the instance type, the presented methods have a function evaluation percentage between $21\%$ and $63\%$ (of the maximum number of function evaluations), and consistently outperform several benchmark methods including the NOMAD solver.

In practical applications, cost functions may be `near-convex', as demonstrated in two application stemming from radiation therapy. In those cases, 1-Opt can provide a near-optimal solution while substantially outperforming the sandwich methods in number of function evaluations. The numerical results indicate that both the choice of bounds (with or without convexity assumption) and the starting point can influence the resulting objective value and the required number of function evaluations. For implementation of the presented solution methods in particular applications, both of these topics require further research.

\section*{Acknowledgements}
We thank Dick den Hertog (University of Amsterdam) for comments that significantly improved this paper, and Thomas Bortfeld (Massachusetts General Hospital and Harvard Medical School) for providing the liver cancer patient data used in \Cref{sec: RT}. We additionally thank Marleen Balvert (Tilburg University), David Craft (Massachusetts General Hospital and Harvard Medical School), Zolt\'{a}n Perk\'{o} (Delft University of Technology), Dick den Hertog and Thomas Bortfeld for inspiring discussions that led to the inception of this paper.

% References
\begingroup
\setstretch{0.9}
\bibliographystyle{apalike}
\small
\bibliography{References}
\endgroup

%\section*{Appendix}
\appendix
\renewcommand{\theequation}{\thesection.\arabic{equation}}
\renewcommand{\thefigure}{\thesection.\arabic{figure}}
\renewcommand{\thetable}{\thesection.\arabic{table}}

\setcounter{equation}{0}
\setcounter{figure}{0}
\setcounter{table}{0}
\sectionfont{\normalsize}
\subsectionfont{\small}

\section{ILP representation} \label{app: ILP}
Resource allocation problem \eqref{eq: knapsack} is equivalent to the following integer (binary) linear programming (ILP) problem:
\begin{subequations} \label{eq: ILP}
\begin{align}
P(F)~=~~~\min_{\bm{y}}~&~ \sum_{i=1}^n \sum_{k=0}^{b_i} F_{i,k}y_{i,k} \\
\text{s.t.}~&~ \sum_{i=1}^n \sum_{k=0}^{b_i} k y_{i,k} = B \label{eq: ILP-total-budget}\\ 
           ~&~ y_{i,k} \in \{0,1\},~\forall k=0,\dotsc,b_i, \forall i=1,\dotsc,n,
\end{align}
\end{subequations}
with $F_{i,k} = f_i(k)$, for all $k=1,\dotsc,b_i$, and $i=1,\dotsc,n$. For known objective coefficients, problem \eqref{eq: ILP} can be solving using any standard ILP solver. The optimal allocation $\bm{x}^{\ast}$ can be recovered from an optimal solution $\bm{y}^{\ast}$ via $\bm{x}_i = \sum_{k=0}^{b_i} k y_{i,k}$ for all $i$.

\section{Radiation therapy data} \label{app: RT-data}
\begin{figure}[htb!]
\hspace*{0.6cm}
\begin{subfigure}[t]{0.48\textwidth}
\hspace*{-0.9cm}
\includegraphics[scale=1]{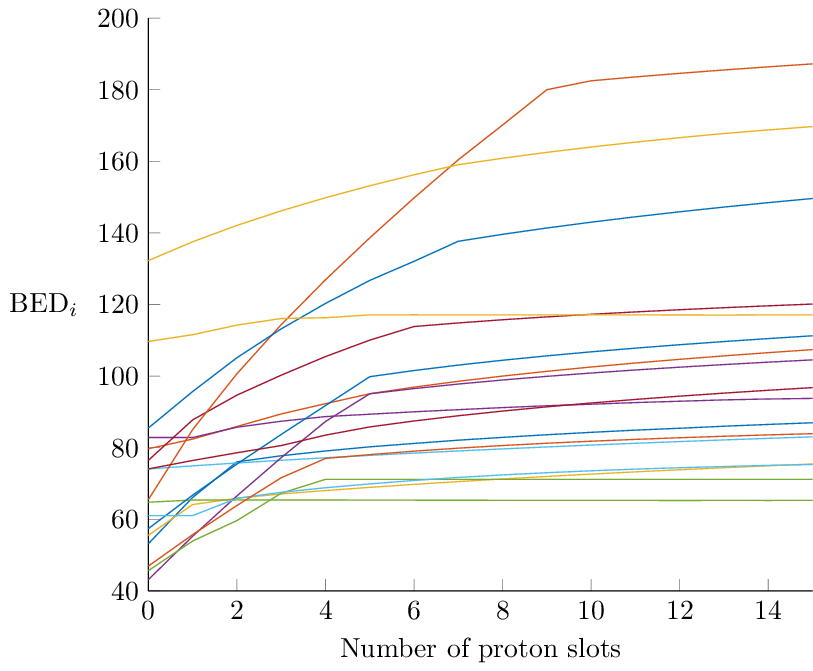}
\caption{\small Payoff functions $\text{BED}_i$ for patient \\$i=1,\dotsc,17$. \label{fig: proton-slot-payoff}}
\end{subfigure}
\hspace*{0.3cm}
\begin{subfigure}[t]{0.48\textwidth}
\hspace*{-0.1cm}
\includegraphics[scale=1]{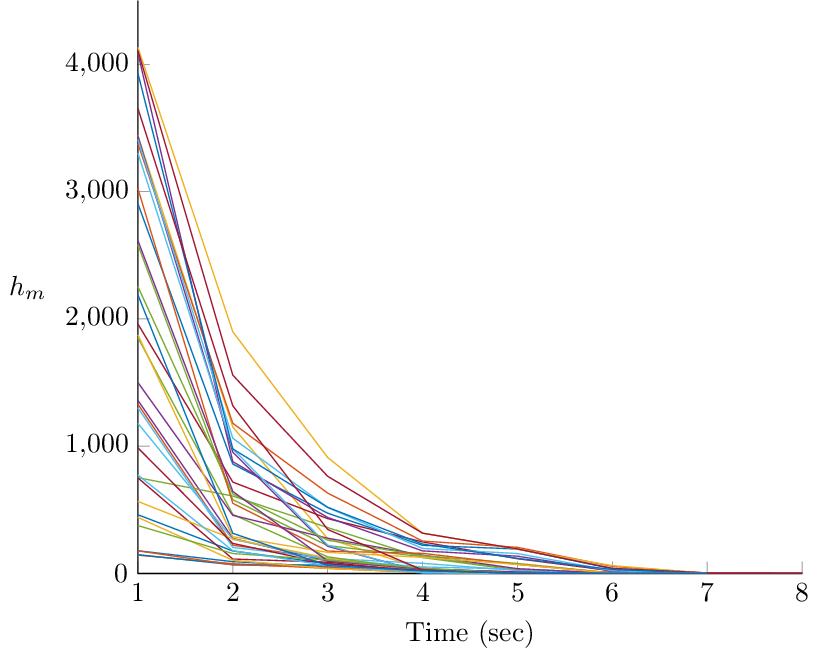}
\caption{\small Fluence map matching inaccuracy functions $h_m$ for arc segments $m=1,\dotsc,36$. \label{fig: VMAT_data}}
\end{subfigure}
\caption{\small Cost functions for the radiation therapy examples in \Cref{sec: RT}.}
\end{figure}

\Cref{fig: proton-slot-payoff} plots the BED for all patients $i=1,\dotsc,17$ as a function of the number of allocated proton slots, for the proton therapy example in \Cref{sec: RT}. \Cref{fig: VMAT_data} plots the fluence map matching inaccuracy $h_m$ for all arc segments $m=1,\dotsc,36$ as a function of time (sec), for the VMAT example in \Cref{sec: RT}.

\setcounter{equation}{0}
\setcounter{figure}{0}
\setcounter{table}{0}
\section{Benchmark methods}\label{app: benchmark}
\subsection{Myopic method} \label{app: myopic}
The myopic method (MY), or greedy method, is an iterative procedure where in each iteration a single item is allocated to a player with the largest immediate gain. For an allocation $\bm{x}^t$ at the start of iteration $t$, the immediate gain of getting an extra item for player $i$ is given by 
\begin{align}
G_i(x_i^t) = f_i(x_i^t) - f_i(x_i^t+1),
\end{align}
which is non-negative, because $f_i$ is non-increasing. The method initializes by computing the immediate gain at zero items for each player $i$, i.e., evaluating $f_i(0)$ and $f_i(1)$. Then, the player $i$ with the highest immediate gain\footnote{In case of ties, the player with the lowest index is chosen.\label{tiebreaker-benchmark}} is assigned the item and its new immediate gain is computed, unless the player $i$ has reached its individual allocation budget $b$. The myopic method starts with $2n$ function evaluations. After that, it performs $B$ allocation iterations. Except for the last iteration, each of these requires at most one additional function evaluation. Thus, the total number of black-box function evaluations for the myopic method is at most $2n+B-1$. 

The total allocation budget $B$ has a large influence on the performance. For instances with a high total allocation budget, we can start with $x_i^0=b_i$ for all players $i$ (i.e., assigning $\bm{b}^{\top}\bm{e}$ items), and remove items until feasible. In particular, in each iteration we remove the item for the player $i$ whose gain for their currently last item is smallest. It is easy to show that if $2B > \bm{b}^{\top}\bm{e}$ this requires less iterations in the worst-case than starting with $x_i^0 = 0$ for all $i$. Pseudocode for the case $2B \leq \bm{b}^{\top}\bm{e}$ is given in \Cref{alg: myopic}. The case with $2B >\bm{b}^{\top}\bm{e}$ is comparable.
\begin{algorithm}[htb!]
\small
\Begin{
Set $t=0$, and set $S = \{1,\dotsc,n\}$, $J=S$\;
Set $x_i^0=0$ for all $i \in S$\;
Evaluate $f_i(0)$ and $f_i(1)$ for all $i \in S$\;
\While{$t < B$}{
Determine $G_i(x_i^t)$ for all $i \in J$\;
	Let $j_t \in \argmax_{i \in J} G_i(x_i^t)$\;
	Set $x_{j_t}^{t+1} = x_{j_t}^{t}+1$\;
	Set $x_i^{t+1} = x_i^t$ for all $i \in S \backslash \{j_t\} $\;
	\eIf{$x_{j_t}^{t}+1 < b_i~ \land~ t<B-1$}{
		Evaluate $f_{j_t}(x_{j_t}^{t}+2)$\;
	}{
		Set $J \leftarrow J \backslash \{ j_t \}$\;
	}
	Set $t \leftarrow t + 1$\;
}
Set $\bm{x}^{\text{ma}} = \bm{x}^t$.
}
\caption{Myopic method for $2B \leq \bm{b}^{\top}\bm{e}$ \label{alg: myopic}}
\end{algorithm}

For convex cost functions the myopic method is exact, as proved in \Cref{lemma: myopic-convexity}.
\begin{lemma} \label{lemma: myopic-convexity}
Let cost functions $f_i: \{0,\dotsc,b_i\} \rightarrow [0,M]$ be convex and non-increasing for all $i=1,\dotsc,n$. Then the myopic solution $\bm{x}^{\text{ma}}$ is optimal to \eqref{eq: knapsack}.
\end{lemma}
\begin{proof}
\small The proof is given for the case $2B \leq \bm{b}^{\top}\bm{e}$, the case $2B > \bm{b}^{\top}\bm{e}$ is similar. Let $G^t$ denote the objective value improvement in iteration $t$ and let $J_t = \{i : x_i^t <b_i \}$, i.e., the set of players eligible for allocation in iteration $t$. For all $t\geq 0$ it holds that
\begin{align}\label{eq: myopic-inequality}
G^t &:= \max_{i \in J_t} G_i(x_i^t) \geq \max_{i \in J_{t+1}} G_i(x_i^{t+1}) = G^{t+1},
\end{align}
where the inequality holds due to convexity of $f_i$ for all $i$ and the fact that $J_{t+1} \subseteq J_t$. With $\bm{x}^{\text{ma}}$ the allocation resulting from the myopic method and $z(\bm{x}^{\text{ma}},f)$ the corresponding objective value, it holds that
\begin{align} \label{eq: myopic-objective}
z(\bm{x}^{\text{ma}},f) =& \sum_{i=1}^n f_i(x_i^{\text{ma}}) = \sum_{i=1}^n f_i(0) - \sum_{t=0}^{B-1} G^t,
\end{align}
i.e., the total cost is the sum of individual player costs at the zero allocation minus the gain in each iteration. From \eqref{eq: myopic-inequality} and \eqref{eq: myopic-objective} it follows that allocating a number of $c$ ($\leq B$) extra items to a subset of players in any way decreases (i.e., improves) the objective value by at most $cG^B$. To maintain feasibility, $c$ items must be removed from the other players, which increases (i.e., deteriorates) the objective value by at least $G^{B-1} + \dotsc + G^{B-c}$, which is larger than $cG^B$. Thus, the current allocation $\bm{x}^{\text{ma}}$ is optimal. \normalsize
\end{proof}

\subsection{Prescient method} \label{app: prescient}
For non-convex cost functions, the myopic method may have poor performance. This poor performance exhibits particularly in those cases where the immediate gain for a player $i$ at $x_i$ is small, but larger gains are possible for higher values of $x_i$. To remedy this, we propose a new heuristic that bases decisions both on the immediate gain and the average gain over the remaining horizon for that player. The prescient method (PR) has the same structure as the myopic method and also allocates items one-by-one. With an allocation $\bm{x}^t$ at the start of iteration $t$, the maximum number of items that can be allocated to player $i$ is
\begin{align}
\beta_i^t(x_i^t) := \min\{b_i, x_i^t+B-t\}.
\end{align}
The average improvement for player $i$ over the interval $[x_i^t,\beta_i^t(x_i^t)]$ is
\begin{subequations}
\begin{align} \label{eq: prescient-average}
A_i(x_i^t) &= \frac{f_i(x_i^t) - u_i(\beta_i^t(x_i^t))}{\beta_i^t(x_i^t) - x_i^t},
\end{align}
\end{subequations}
where the upper bound value $u_i(\beta_i^t(x_i^t))$ is used because $f_i(\beta_i^t(x_i^t))$ need not be evaluated. Thus, the reported average gain is conservative. The score of player $i$ at $x_i^t$ is defined as
\begin{align} \label{eq: prescient-score}
s_i(x_i^t) := \max \{ G_i(x_i^t), A_i(x_i^t) \}.
\end{align}
In each iteration, the item is allocated to the player $i$ with the currently highest score\footref{tiebreaker-benchmark}. The prescient method starts with $3n$ function evaluations. After that, it performs $B$ allocation iterations. Except for the last iteration, each of these requires at most one additional function evaluation. Thus, the total number of black-box function evaluations for the prescient method is at most $3n+B-1$.

Similar to the myopic method, for instances with a high total allocation budget we can start with $x_i^0=b_i$ for all players $i$ (i.e., assigning $\bm{b}^{\top}\bm{e}$ items), and remove items until feasible. Similar to the myopic method, removing items requires fewer iterations (in the worst-case) than adding items if $2B \leq \bm{b}^{\top}\bm{e}$. Pseudocode for the case $2B \leq \bm{b}^{\top}\bm{e}$ is given in \Cref{alg: prescient}.
\begin{algorithm}[h]
\small
\Begin{
Set $t=0$, and set $S = \{1,\dotsc,n\}$, $J=S$\;
Set $x_i^0 = 0$ for all $i \in S$\;
Evaluate $f_i(0)$, $f_i(1)$ and $f_i(b)$ for all $i \in S$\;
\While{$t < B$}{
	Determine $s_{i}(x_i^t)$ for all $i\in J$\;
	Let $j_t \in \argmax_{i \in J} s_i(x_i^t)$\;
	Set $x_{j_t}^{t+1} = x_{j_t}^{t}+1$\;
	Set $x_i^{t+1} = x_i^t$ for all $i \in S \backslash \{j_t\} $\;
	\uIf{$x_{j_t}^{t}+1 < b ~\land~ t < B-1$}{
			Evaluate $f_{j_t}(x_{j_t}^{t}+2)$\;
	}
	\uElseIf{$x_{j_t}^{t}+1 = b$}{	
		Set $J \leftarrow J \backslash \{ j_t \}$\;
	}
	Set $t \leftarrow t + 1$\;
}
Set $\bm{x}^{\text{pa}} = \bm{x}^t$\;
}
\caption{Prescient method for $2B \leq \bm{b}^{\top}\bm{e}$ \label{alg: prescient}}
\end{algorithm}

For convex cost functions the prescient method is exact, as proved in \Cref{lemma: prescient-convexity}.
\begin{lemma} \label{lemma: prescient-convexity}
Let cost functions $f_i: \{0,\dotsc,b_i\} \rightarrow [0,M]$ be convex and non-increasing for all $i=1,\dotsc,n$. Then the prescient solution $\bm{x}^{\text{pa}}$ is optimal to \eqref{eq: knapsack}.
\end{lemma}
\begin{proof}
\small
Due to convexity of cost functions $f_i$, it holds that $G_i(x_i) \geq A_i(x_i)$ and, consequently, $s_i(x_i) = G_i(x_i)$ for all feasible $x_i$ and all $i=1,\dotsc,n$. Thus, the prescient method uses the same allocation rule as the myopic method in each iteration, and optimality follows from \Cref{lemma: myopic-convexity}. \normalsize
\end{proof}
Both presented benchmark methods, the myopic and the prescient method, are optimal for convex cost functions, and the latter requires an equal or higher number of function evaluations than the former. However, for non-convex cost functions the prescient method is expected to perform better because it also accounts for possible non-convexities via the score function~\eqref{eq: prescient-score}. 

In the numerical experiments, we use \Cref{alg: myopic,alg: prescient} if $2B \leq \bm{b}^{\top}\bm{e}$. Otherwise, we use their counterparts starting with each player $i$ allocated $b_i$ items, and use the described procedure to remove items instead.

Lastly, we note that simple adaptations may improve the myopic and prescient method. For example, the direct gain of a player need not be computed if the (upper bound on the) total possible gain for the player is lower than the direct gain of another player. For ease of exposition, we do not incorporate such adaptations.

\end{document}